\documentclass[a4paper,american]{amsart}

\usepackage{times,latexsym,amssymb}
\usepackage{amsmath,amsthm,bm}
\usepackage{color}
\usepackage[makeroom]{cancel}
\usepackage{mathrsfs}

\usepackage{amssymb}
\usepackage{dsfont}
\usepackage{stmaryrd}
\usepackage{esint}
\usepackage{changes}

\usepackage{xcolor}
\usepackage[colorlinks,pdfpagelabels,pdfstartview = FitH,bookmarksopen = true,bookmarksnumbered = true,linkcolor = blue,plainpages = false,hypertexnames = false,citecolor = red,pagebackref=false]{hyperref}
\usepackage{enumitem}
\usepackage{changes}
\usepackage{comment}
\allowdisplaybreaks

\newtheorem{theorem}{Theorem}
\newtheorem{lemma}[theorem]{Lemma}
\newtheorem{proposition}[theorem]{Proposition}

\theoremstyle{definition}
\newtheorem{definition}[theorem]{Definition}

\theoremstyle{remark}
\newtheorem{remark}[theorem]{Remark}

\numberwithin{theorem}{section}
\numberwithin{equation}{section}

\newcommand{\N}{\mathbb{N}}
\newcommand{\R}{\mathbb{R}}

\renewcommand{\d}{\mathrm{d}}
                  
\newcommand{\dx}{\mathrm{d}x}
\newcommand{\dy}{\mathrm{d}y}

\newcommand{\dt}{\mathrm{d}t}

\newcommand{\loc}{\mathrm{loc}}

\renewcommand{\epsilon}{\varepsilon}
\renewcommand{\rho}{\varrho}

\newcommand{\dsty}{\displaystyle}
\newcommand{\al}{\alpha}

\newcommand{\gm}{\gamma}

\newcommand{\Dl}{\Delta}

\newcommand{\Lm}{\Lambda}

\newcommand{\om}{\omega}
\newcommand{\Om}{\Omega}

\newcommand{\z}{\zeta}

\newcommand{\bth }{\boldsymbol{\tau}_h}

\def\YYint#1#2#3{{\setbox0=\hbox{$#1{#2#3}{\iint}$}
    \vcenter{\hbox{$#2#3$}}\kern-.51\wd0}}

\def\Xint#1{\mathchoice
{\XXint\displaystyle\textstyle{#1}}%
{\XXint\textstyle\scriptstyle{#1}}%
{\XXint\scriptstyle\scriptscriptstyle{#1}}%
{\XXint\scriptscriptstyle\scriptscriptstyle{#1}}%
\!\int}
\def\XXint#1#2#3{{\setbox0=\hbox{$#1{#2#3}{\int}$ }
\vcenter{\hbox{$#2#3$ }}\kern-0.555\wd0}}

\def\fint{\Xint-}

\DeclareMathOperator{\spt}{spt}

\DeclareMathOperator*{\osc}{osc}
\DeclareMathOperator*{\esssup}{ess\,sup}
\DeclareMathOperator*{\essinf}{ess\,inf}
\DeclareMathOperator*{\essosc}{ess\,osc}
\DeclareMathOperator{\Div}{div}

\def\Xint#1{\mathchoice
    {\XXint\displaystyle\textstyle{#1}}%
    {\XXint\textstyle\scriptstyle{#1}}%
    {\XXint\scriptstyle\scriptscriptstyle{#1}}%
    {\XXint\scriptscriptstyle\scriptscriptstyle{#1}}%
    \!\int}
\def\XXint#1#2#3{\setbox0=\hbox{$#1{#2#3}{\int}$}
    \vcenter{\hbox{$#2#3$}}\kern-0.5\wd0}

\def\dashint{\Xint{\raise4pt\hbox to7pt{\hrulefill}}}

\def\Xiint#1{\mathchoice
    {\XXiint\displaystyle\textstyle{#1}}%
    {\XXiint\textstyle\scriptstyle{#1}}%
    {\XXiint\scriptstyle\scriptscriptstyle{#1}}%
    {\XXiint\scriptscriptstyle\scriptscriptstyle{#1}}%
    \!\iint}
\def\XXiint#1#2#3{\setbox0=\hbox{$#1{#2#3}{\iint}$}
    \vcenter{\hbox{$#2#3$}}\kern-0.5\wd0}
\def\biint{\Xiint{-\!-}}

\def\XXiiint#1#2#3{\setbox0=\hbox{$#1{#2#3}{\iint}$}
    \vcenter{\hbox{$#2#3$}}\kern-0.5\wd0}


\begin{document}

\title[$C^1$-estimates in nonlocal parabolic problems]{$C^1$-estimates in nonlocal parabolic problems
}
\date{\today}

\author[V. B\"ogelein]{Verena B\"{o}gelein}
\address{Verena B\"ogelein\\
Fachbereich Mathematik, Universit\"at Salzburg\\
Hellbrunner Str. 34, 5020 Salzburg, Austria}
\email{verena.boegelein@plus.ac.at}

\author[N. Liao]{Naian Liao}
\address{Naian Liao\\
Fachbereich Mathematik, Universit\"at Salzburg\\
Hellbrunner Str. 34, 5020 Salzburg, Austria}
\email{naian.liao@plus.ac.at}

\author[C.~Stanko]{Calvin Stanko}
\address{Calvin Stanko\\
Fachbereich Mathematik, Universit\"at Salzburg\\
Hellbrunner Str.~34, 5020 Salzburg, Austria}
\email{calvin.stanko@plus.ac.at}

\begin{abstract}
     Under a minimal tail condition, local weak solutions to a class of nonlinear nonlocal parabolic equations modeled on the fractional heat equation possess continuous spatial derivatives.
\end{abstract}

\maketitle

\tableofcontents

\section{Introduction}
    A special case of our main result is as follows. 
    Let $s\in(0,1)$ and $u$ be a local weak solution to the fractional heat equation
    \begin{equation*}
        \partial_t u +(-\Dl )^s u=0 \quad\text{in}\quad Q_{4}= B_{4}\times(-4^{2s},0].
    \end{equation*}
    Then, its spatial derivatives $D_i u$, $i \in \{1,\dots,N\}$ are locally continuous in $Q_{4}$. Moreover, we have the local boundedness estimate:  
    \begin{align}\label{sup-est-sample}
        \sup_{Q_1} |D u|  \lesssim  \int^0_{-2^{2s}}\int_{\R^N}\frac{|u(x,t)|}{(1+|x|)^{N+2s}}\,\dx\dt,
    \end{align}
    and the modulus of continuity at each point $(x_o,t_o)\in Q_1$:
    \begin{align}\label{osc-est-sample}
        \osc_{Q_{r}(x_o,t_o)} D_{i} u \lesssim \boldsymbol{\om} \, r^\beta +  \int_{t_o-r^s}^{t_o}\int_{\R^N 
        }\frac{|u(x,t)|}{(1 + |x-x_o|)^{N+2s}}\,\dx\dt  
    \end{align}
    for any $r\in (0,1)$,  where   $\beta \in(0,1)$ depends on $\{s, N\}$ and
    \begin{align*}
        \boldsymbol{\om}= \int_{t_o-2^{2s}}^{t_o}\int_{\R^N}\frac{|u(x,t)|}{(1+|x-x_o|)^{N+2s}}\,\dx\dt.
    \end{align*}

A distinctive feature of nonlocal equations is that the local behavior of solutions is influenced by their values at arbitrarily large distances. Consequently, any local regularity estimate must necessarily involve suitable global quantities measuring the long-range interactions encoded by the operator. For the fractional heat equation, the weighted space-time integral appearing on the right-hand side of \eqref{sup-est-sample}--\eqref{osc-est-sample} is precisely the natural quantity entering the definition of local weak solutions. The main point of the present paper is that this quantity  completely controls their interior regularity up to $C^1$-level.

 Kassmann and Weidner \cite{KW_24}  proved the local boundedness estimate  of the form~\eqref{sup-est-sample} for {\it solutions} to a general equation modeled on the fractional heat equation, having measurable coefficients; see also \cite{CLD_2016} in a different context. They also established Hölder continuity of {\it solutions} under the stronger assumption of $L^{1+\varepsilon}$-time-integrability and showed that Hölder continuity may fail under the natural $L^1$-time-integrability assumption. This left open the question of whether continuity itself still follows from the minimal integrability required in the weak formulation. The question was answered affirmatively in \cite{Liao_24_3}, where a quantitative modulus of continuity of the form \eqref{osc-est-sample} was established for {\it solutions}.

The present paper shows that the same natural global quantity governs the regularity of the {\it gradient of solutions}. More precisely, under the minimal $L^1$-time-integrability assumption, we prove that the gradient exists, is continuous, and satisfies quantitative local estimates depending only on the corresponding tail of the solution. Thus, no additional global assumptions are required to obtain   interior regularity up to $C^1$-level. Beyond the linear fractional heat equation, the arguments developed here apply to a class of nonlinear nonlocal parabolic equations, for which they provide a unified gradient regularity theory.

 \subsection{Main results}
 Let $\Omega \subset \mathbb{R}^{N}$, $N \ge 2$ be an open and bounded set; for $T > 0$ define the space-time cylinder $\Omega_T := \Omega \times (0, T]$. 
Consider the nonlinear nonlocal parabolic equation
\begin{align}\label{eq:PDE}
\partial_t u(x,t)+\mathscr{L} u(x,t)= 0 \quad\text{in}\>\> \Omega_T,
\end{align}
where  
\begin{align*}
    \mathscr{L} u(x,t) := \text{ p.v.} \int_{\mathbb{R}^{N}} \Phi \bigg( \frac{u(x,t) - u(y,t)}{| x-y |^{s}} \bigg) \frac{\dy}{\vert x-y\vert^{N+s}}
\end{align*}
for some $s \in (0,1)$, and the nonlinearity $\Phi \colon \mathbb{R} \to \mathbb{R}$ is an odd function such that for any $\xi_1,\, \xi_2 \in \mathbb{R}$ and some $\Lambda \ge 1$, we have the growth and coercivity conditions:
\begin{equation}\label{assumption1:nonlinear_integrand_kernel}
\begin{aligned}
    | \Phi (\xi_1) - \Phi (\xi_2) | &\le \Lambda | \xi_1 - \xi_2 | \\
     ( \Phi (\xi_1) - \Phi (\xi_2)) ( \xi_1 - \xi_2 ) &\ge \Lambda^{-1} | \xi_1 - \xi_2 |^2.
\end{aligned}
\end{equation}
Our main result is as follows.
\begin{theorem}\label{Thm:main}
        Let $s\in(0,1)$ and let $u$ be a local weak solution to \eqref{eq:PDE}--\eqref{assumption1:nonlinear_integrand_kernel} in the sense of Definition~\ref{def:weak_solution}.
    Then, its spatial derivatives $D_i u$, $i \in \{1,\dots,N\}$ are locally continuous functions in $\Om_T$. Moreover,  there exist $C>1$ and $\beta \in(0,1)$ depending on $\{ s, N, \Lm\}$, such that for any $Q_{4R}(y,\tau)\subset \Om_T$, 
    we have the local boundedness estimate:
    \begin{align*}
        \sup_{Q_R(y,\tau)} |D u| \le \frac{C}{R} \int^{\tau}_{\tau-(2R)^{2s}}\int_{\R^N}\frac{|u(x,t)|}{(R+|x|)^{N+2s}}\,\dx\dt,
    \end{align*}
    and the modulus of continuity at each point $(x_o,t_o)\in Q_R(y,\tau)$:
    \begin{align*}
        \osc_{Q_{ r}(x_o,t_o) } D_{i} u \le  C \bigg[  \boldsymbol{\om} \Big( \frac{r}{R} \Big)^{\beta} + \frac{1}{R} \int_{t_o - (rR)^{s}}^{t_o} \int_{\mathbb{R}^{N} } \frac{| u(x,t) |}{(R+| x - x_o |)^{N+2s}} \,\dx\dt \bigg]
    \end{align*}
    for any $r\in (0,R)$, where
    \begin{align}\label{def:final_omega}
        \boldsymbol{\om}= \frac{1}{R}\int_{t_o-(4R)^{2s}}^{t_o}\int_{\R^N}\frac{|u(x,t)|}{(R+|x-x_o|)^{N+2s}}\,\dx\dt.
    \end{align}
\end{theorem}

\begin{remark}
    This result can be viewed as a nonlocal analogue of the well-known $C^{1,\alpha}$-regularity of local weak solutions to the quasilinear parabolic equation
\begin{align*}
    \partial_t u - \Div \mathbf{A}(D u)=0
\end{align*}
where the $C^1$-vector field $\mathbf{A}:\R^N\to\R^N$ satisfies proper coercivity and growth conditions;
cf.~\cite[Chapter~V, \S\,3, Theorem 3.1]{LSU}.
\end{remark}

\begin{remark}
The obtained $C^1$-regularity is sort of ``end-point'' regularity.
    Under a more stringent $L^{1+\epsilon}$-time-integrability of the tail, Theorem~\ref{Thm:main} yields $C^{1,\al}$-regularity, and hence recovers \cite[Theorem~1.6]{DKLN_2025}. Our study also gives further improvement, which will appear elsewhere.
\end{remark}

\begin{remark}
    Nonlocal higher regularity in the time variable can be studied under additional conditions. However, under the minimal tail condition, one cannot expect further regularity  in the time variable beyond absolute continuity.
\end{remark}

\subsection{The strategy}
The method of difference quotient, a traditional barehanded technique in the analysis of partial differential equations, has found new applications in various nonlocal problems. See~\cite{BDLGMS24, BDLGM24-2, BL-17, BLS-18, CCV_2011, DKLN_ann, DKLN_2025, Diening_Kim_Lee_Nowak_24} for instance.

In a previous work \cite{DKLN_2025}, the authors, within the same framework as this paper, assumed additional time-integrability in the tail condition and thereby established $C^{1,\alpha}$-regularity. Roughly speaking, they considered the difference quotient $\bth u/|h|^{\gm}$, $\gm\in(0,1]$ and leveraged available Hölder estimates for linear equations with measurable kernels, which hold only when the additional tail condition is given, thereby inductively increasing the Hölder exponent of $u$. The method is built on finite iterations in Hölder spaces. Therefore, we cannot  adapt their scheme to attack the ``end-point'' gradient regularity under the natural tail condition. 

To reach our goal, we must insist on Sobolev spaces, examine the equation satisfied by $\bth u/|h|$, and perform various more refined iterative arguments within a framework capable of handling the natural tail. In short, our regularity-upgrade-scheme can be illustrated by
\begin{align*}
    L^2_t W^{s,2}_x \cap L^1_tL^1_{2s,x}\longrightarrow L^2_t W^{1,2}_x \cap L^1_tL^1_{2s,x} \longrightarrow L^\infty_t W^{1,\infty}_x \cap L^1_tL^1_{2s,x}.
\end{align*}
In the first step, only the fractional differentiability is improved while the integrability in both time and space remains the same. This step adapts some arguments from \cite{BDLGMS24,BDLGM24-2}. In the second step, we use De Giorgi's iteration to show that $\|\bth u/|h|\|_{L^\infty_t L^{\infty}_x}$ can be bounded by $\|\bth u/|h|\|_{L^2_t L^{2}_x}$ and the natural tail of $u$, thereby establishing the spatial Lipschitz regularity.

Continuity of the gradient entails further delicate analysis. Again, we examine the equation satisfied by $\bth u/|h|$, use some techniques developed in \cite{Liao_24_2, Liao_24_3}, and derive a modulus of continuity of $\bth u/|h|$ that is uniform in $h$. Passing to the limit, we eventually achieve the desired modulus of continuity for $Du$.

Last but not least, we would like to point out that even in the elliptic case, our route of upgrading the regularity from $W^{s,2}$ to $C^{1,\al}$ is new, comparing to \cite{Diening_Kim_Lee_Nowak_24} within the same framework. The latter again features finite iterations in Hölder spaces. Moreover, we do not use so-called “localization” techniques, but rather handle the tail terms directly.

\

\subsection*{Acknowledgments}  
 This research was set out in May, 2025. No AI tools were used in the process.
It was funded in whole or in part by the Austrian Science Fund (FWF) [10.55776/PAT1850524]. For open access purposes, the author has applied a CC BY public copyright license to any author accepted manuscript version arising from this submission.

\section{Notation and definitions}

Our notation is standard in the field. We denote the ball in $\R^N$ with radius $R$ and center $x_o$ by $B_R(x_o)$; a space-time cylinder at $z_o=(x_o,t_o)$ is
\begin{align*}
    Q_{R,S}(z_o):=B_R(x_o)\times (t_o-S,t_o]  
\end{align*}
and from time to time, we also use
\begin{align*}
    \mathbb{Q}_{R,S}(z_o):=B_R(x_o)\times B_R(x_o)\times (t_o-S,t_o].
\end{align*}
The parabolic cylinders are   $Q_R(z_o):=Q_{R,R^{2s}}(z_o)$ and $\mathbb{Q}_{R}(z_o):=\mathbb{Q}_{R,R^{2s}}(z_o)$.

For $s \in (0,1)$ define the weighted Lebesgue space
\begin{align*}
    L^{1}_{2s} ( \mathbb{R}^{N}):=
    \bigg\{ u\in L^1_{\loc}(\R^N) \,:\,
    \int_{\R^N} \frac{|u(x)|}{1+|x|^{N+2s}} \,\dx < \infty\bigg\}.
\end{align*}
The global behavior of a function belonging to the space 
$L_{2s}^{1}(\mathbb{R}^{N})$
can be measured by the tail, which is defined by
\begin{align}\label{def:elliptic_tail}
    \mathrm{Tail}(u(t); B_R(x_o)) := R^{2s} \int_{\mathbb{R}^{N} \setminus B_{R}(x_o) } \frac{\vert u(x,t) \vert}{\vert x - x_o \vert^{N+2s}} \,\dx .
\end{align}
When dealing with a time-dependent function $u\in L^1(t_o-S, t_o;L^{1}_{2s} ( \mathbb{R}^{N}))$ we also use
\begin{align*}
    \mathrm{Tail}(u; Q_{R,S}(z_o)) := \int_{t_o-S}^{t_o} \int_{\mathbb{R}^{N} \setminus B_{R}(x_o) } \frac{\vert u(x,t) \vert}{\vert x - x_o \vert^{N+2s}} \,\dx\dt .
\end{align*}
Moreover, 
we introduce the abbreviations
\begin{align}\label{def:D-tail-ellip}
    \mathcal{T} (t;B_R(x_o)) := \Vert  u(t) \Vert_{L^{\infty}(B_{R}(x_o))} + \mathrm{Tail} (u(t); B_R(x_o) )
\end{align}
and
\begin{align}\label{def:D-tail-par}
    \mathcal{T} ( Q_{R,S}(z_o)) :=  \Vert  u \Vert_{L^{\infty}(Q_{R,S}(z_o))} + \mathrm{Tail} (u; Q_{R,S}(z_o) ).
\end{align}
In the case the center is the origin or is clear from the context, we omit it. 

Throughout this paper, we let
\begin{align*}
    \mathfrak{data}:=\{s, N, \Lm\},
\end{align*}
and use $C = C(\mathfrak{data})$ as a generic positive constant, which can be determined by the $\mathfrak{data}$ only and can vary from one estimate to another.

We end this section by defining local weak solutions to~\eqref{eq:PDE}--\eqref{assumption1:nonlinear_integrand_kernel} as follows.

\begin{definition}\label{def:weak_solution}(Local weak solution)
Let  $s \in (0,1)$. A function 
\begin{align*}
        u \in C_{\loc} (0,T; L_{\loc}^{2} ( \Omega ) ) \cap L_{\loc}^{2} (0,T; W_{\loc}^{s,2}(\Omega)) \cap 
        L^1_\loc (0,T;L^1_{2s}(\R^N))
\end{align*}
is a local, weak solution to~\eqref{eq:PDE}--\eqref{assumption1:nonlinear_integrand_kernel}, 
if for every sub-interval $[t_1, t_2] \subset (0,T]$, we have
\begin{align}\label{eq:weak_formulation}
    0=&\int_{\Omega} \phi\, u \,\dx \bigg|_{t_1}^{t_2}  
    -
    \int_{t_1}^{t_2} \int_{\Omega} \partial_t \phi\, u\,\dx \dt \nonumber\\
    &+\int_{t_1}^{t_2} \int_{\mathbb{R}^N}\int_{\mathbb{R}^N} \Phi \bigg( \frac{u(x,t)-u(y,t)}{\vert x-y \vert^{s}} \bigg) \frac{\phi(x,t)-\phi(y,t)}{\vert x-y \vert^{N+s}} \,\dy \dx  \dt ,
\end{align}
for any 
$\phi\in L_{\loc}^2(0,T;W_{0}^{s,2}(\Omega))\cap W_{\loc}^{1,2}(0,T;L^2(\Omega))$ such that $\spt\phi(\cdot, t)\subset\Omega$. 
\end{definition}

\section{Preliminaries}

\subsection{Algebraic inequalities and auxiliary lemmata}

At several stages we will rely on algebraic inequalities.

\begin{lemma}\label{lem:ineq_algebraic_a_b_c_d_e_f_p_euqal_to_2}
    Let $a,b,e,f \in \mathbb{R}$.
    Then we have
    \begin{align*}
        (a-b)\big(ae^2-bf^2\big)
        \ge 
        \tfrac14 (a-b)^2\big(e^2+f^2\big) -
        \tfrac12 (a+b)^{2} \vert e - f \vert^{2}.
    \end{align*}
\end{lemma}

\begin{proof}
The left-hand side can be rewritten as
\begin{align*}
    (a-b)\big(ae^2-bf^2\big)
    &=
    \tfrac12 (a-b)^2\big(e^2+f^2\big) + 
    \tfrac12 (a-b)(a+b)\big(e^2-f^2\big) .
\end{align*}
The second term on the right-hand side can be bounded as
\begin{align*}
    \tfrac12 \big|(a-b)(a+b)\big(e^2-f^2\big)\big|
    &\le 
    \tfrac18(a-b)^2(e+f)^2 + \tfrac12(a+b)^2(e-f)^2 \\
    &\le 
    \tfrac14(a-b)^2\big(e^2+f^2\big) + \tfrac12(a+b)^2(e-f)^2.
\end{align*}
Inserting it above yields the claimed inequality.
\end{proof}

\subsection{An integral estimate}

Next, we will introduce an estimate for 
the tail introduced as in~\eqref{def:elliptic_tail}.
This result can be found in~\cite[Lemma 2.3]{BLS-18} and will be later applied on each time-slice of our space-time cylinder.

\begin{lemma}\label{lem:prelim_tail_estimate}
    For $s \in (0,1)$, $B_{R} \equiv B_{R}(x_o) \Subset \Omega$ and $r \in (0,R)$, we have whenever $u \in L_{\loc}^{\infty}(\Omega) \cap L_{2s}^{1}(\mathbb{R}^{N})$ that  
    \begin{align*}
        \mathrm{Tail}(u; B_r  ) \le C(N) \Big( \frac{R}{r} \Big)^{N} \Big[ \mathrm{Tail}(u; B_R ) + \Vert u \Vert_{L^{\infty}(B_{R})} \Big].
    \end{align*}
\end{lemma}

\subsection{Finite differences and fractional Sobolev spaces}

Suppose that $\Omega \subset \mathbb{R}^{N}$ is open and $h \in \mathbb{R}^{N}$ is a direction vector. Then, let $\Omega_{h}:= \{ x \in \Omega: x+h \in \Omega \}$. This enables us to define for a measurable function $v\colon \Omega \to \mathbb{R}$ the finite difference operator $\bth\colon L^{1}(\Omega) \to L^{1}(\Omega_h)$ as
\begin{align*}
    \bth  v(x,t):=v(x+h,t) - v(x,t).
\end{align*}
Instead of $v(x+h,t)$ we will often use the abbreviation $v_{h}(x,t):=v(x+h,t)$. For a fixed unit vector $e \in \mathbb{R}^{N}$ and $h\in\R$, we introduce the notation
\begin{align*}
    \bth^{(e)}  v(x,t):=v(x+he,t) - v(x,t)
\end{align*}
and write $\bth^{(i)}  v$ in the case of a canonical basis vector $e_{i} \in \mathbb{R}^{N}$.
 The following property of finite differences are well known, cf.~\cite[Proposition~20.3]{DB-RA}. For the rest of the section, we will omit the center of balls.

\begin{lemma}\label{lem:bound_of_finite_difference_implies_bound_of_first_derivative}
    Assume that $M>0$. Then any $v \in L^{2}(B_R)$ satisfying 
    \begin{align}\label{ineq:condition_that_finite_difference_is_bounded_by_h}
        \int_{B_{\frac12 R}} \vert \bth^{(i)} v \vert^{2} \,\dx \le M^2 \vert h \vert^{2}, 
    \end{align}
    for all $0 < \vert h \vert \le \frac12 R$, is weakly differentiable with respect to $x_{i}$ on $B_{\frac12 R}$ and
    \begin{align*}
        \int_{B_{\frac12 R}} \vert D_{i} v \vert^{2} \,\dx \le M^{2}.
    \end{align*}
    If~\eqref{ineq:condition_that_finite_difference_is_bounded_by_h} holds for any $e_i$, then $v \in W^{1,2}(B_{\frac12 R})$.
\end{lemma}

The next result displays an estimate for finite differences of $W^{\gamma,2}$-functions, cf.~\cite[Proposition 2.6]{BL-17}. 

\begin{lemma}\label{lem:ineq_for_finite_differences_of_sobolev_functions}
    For $\gamma \in (0,1)$ 
    there exists a constant $C=C(N)$ such that for $v \in W^{\gamma,2}(B_R)$ and $h\in \mathbb{R}^{N} \setminus \{0\}$ with $\vert h \vert \le \frac12 R$ we have the estimate
    \begin{align*}
        \int_{B_{\frac12 R}} \vert \bth  v \vert^{2} \,\dx \le C \vert h \vert^{2 \gamma } \Big[ [v]_{W^{\gamma,2}(B_R)}^{2} +    R^{-2\gamma }  \Vert v \Vert_{L^{2}(B_R)}^{2}  \Big].
    \end{align*}
\end{lemma}

The following lemma, important in handling second order differences, is essentially contained in \cite[Chapter 5]{Stein}; see also~\cite[Lemma~2.17]{BDLGMS24}, \cite[Proposition 2.4]{BL-17}, as well as~\cite[Theorem 1.1]{Domokos_2004_2}. 

\begin{lemma}\label{lem:domokos}
    Suppose that $\gamma \in (0,1)$, and $M \ge 0$. Assume that $v \in  L^{2}(B_{R})$  
    satisfies
    \begin{align*}
        \int_{B_{\frac12 R}} \vert \bth^2  v  \vert^{2} \,\dx \le M^{2} \vert h \vert^{2 \gamma },
    \end{align*}
    for all $h \in \mathbb{R}^{N} \setminus \{ 0 \}$ satisfying $\vert h \vert \le \frac14 R$. Then, there exist an absolute constant $C$ such that 
    \begin{align*}
           \int_{B_{\frac12 R}} \vert \bth  v  \vert^{2} \,\dx \le C \vert h \vert^{2 \gamma } \bigg[ \frac{M^{2}}{(1-\gm)^2} +  R^{- 2 \gamma }  \int_{B_R} \vert v \vert^{2} \,\dx  \bigg]
    \end{align*}
    for any $h \in \mathbb{R}^{N} \setminus \{ 0 \}$ satisfying $0 < |h| \le \frac18 R$.
\end{lemma}

At this point we will also introduce a tail estimate for finite differences, cf.~\cite[Lemma 3.1]{BDLGMS24}. Moreover, note that due to the nature of the given problem we will apply it on each time-slice and with the choice $p=2$.

\begin{lemma}\label{lem:first_tail_estimate}
    Let $s\in(0,1)$, $R>0$, and $r\in(0,R)$. Consider the function 
    $$
    u \in L^\infty(B_{2R} )\cap L_{2s}^{1}(\mathbb{R}^{N}).
    $$ 
    Then, there exists a constant $C=C(N,s)$ such that for a.e. $x \in B_{\frac12 (R+r)} $ and any $h \in \mathbb{R}^{N} \setminus \{0\}$ satisfying $0<\vert h \vert \le \frac14 (R-r)$, we have
        \begin{align*}
        \bigg|\int_{\R^N\setminus B_R} 
    \frac{\bth u(y) }{|x-y|^{N+2s}} \,\dy \bigg|
    \le \frac{C|h|}{R^{2s+1}}\Big(\frac{R}{R-r}\Big)^{N+2s+1} \mathcal{T},
    \end{align*}
    where
    \begin{align*}
        \mathcal{T}:= \Vert u \Vert_{L^{\infty}(B_{2R} )}+\mathrm{Tail}(u; B_{2R}).
    \end{align*}
\end{lemma}

In the following we state an embedding for fractional Sobolev spaces. 
For the proof we refer to~\cite[Proposition A.3]{Liao_24_2}.

\begin{proposition}\label{prop:parabolic_embedding}
    Let $s \in (0,1)$ and $R>0$. Consider the function
    \begin{align*}
        u \in L^{2} (t_1, t_2; W^{s,2} (B_R) ) \cap L^{\infty}(t_1, t_2; L^{2}(B_R)),
    \end{align*}
    which is supported in $B_{\frac12 R}$  for a.e. $t \in (t_1, t_2)$. Then, 
    \begin{align*}
         &\int_{t_1}^{t_2} \int_{B_R} \vert u\vert^{\frac{2(N+2s)}{N}} \,\dx\dt\\
         &\quad\le C \bigg[  R^{2s} \int_{t_1}^{t_2} \int_{B_R} \int_{B_R} \frac{\vert u(x,t) - u(y,t) \vert^2}{\vert x - y \vert^{N+2s}} \,\dx\dy\dt 
         + 
         \int_{t_1}^{t_2} \int_{B_R} \vert u  \vert^{2} \,\dx\dt \bigg]\\ &\qquad\quad\cdot\bigg[ \esssup\limits_{t \in (t_1, t_2)} \fint_{B_R\times\{t\}} \vert u  \vert^{2} \,\dx\bigg]^{\frac{2s}{N}}
    \end{align*}
    holds true for some $C=C(N,s)$.
\end{proposition}

\section{Energy inequalities}

In this section, we present two kinds of energy estimates regarding the finite difference $\bth u $.
The first one will be used to upgrade the preset $W^{s,2}$-regularity of weak solutions to the $W^{1,2}$-regularity in the space variable.

\begin{lemma}\label{cor:second_energy_inequality}
    Let $s\in(0,1)$. There exists a constant $C=C(\mathfrak{data})>0$ such that, if $u$
is a locally bounded, local weak solution to~\eqref{eq:PDE} in the sense of Definition~\ref{def:weak_solution},  
then for any $0<r<R$ and $0<\tau<S$ satisfying $Q_{2R,S}(z_o)  \Subset \Omega_T$, and any $h \in \mathbb{R}^{N} \setminus \{ 0 \}$ with $0<\vert h \vert \le \frac14 (R-r)$, we have the energy inequality 
\begin{align*}
        \esssup\limits_{t\in (t_o-\tau, t_o)}&\int_{B_{r}\times\{t\}} \vert \bth u \vert^{2}   \,\dx + \iiint_{\mathbb{Q}_{r,\tau}} \frac{\vert \bth u(x,t)  - \bth u(y,t)  \vert^{2} }{\vert x-y \vert^{N+2s}} \,\dx\dy\dt \nonumber\\
        &\le 
        C \bigg[\frac{1}{R^{2s}} \Big( \frac{R}{R-r} \Big)^{N+2s} + \frac{1}{S-\tau}\bigg] 
        \iint_{Q_{R,S}} \vert \bth u\vert^{2} \,\dx\dt \\
        &\phantom{=}+ 
        \frac{C \vert h \vert^2}{R^{4s+2}} \Big( \frac{R}{R-r} \Big)^{2N+4s+2} R^{N} \bigg[ \int_{t_o-S}^{t_o} \mathcal{T}(t;B_{2R}) \,\dt \bigg]^2 ,
\end{align*}
where we omitted $x_o$, $z_o$ and  defined $\mathcal{T} (t;B_{2R})$ as in \eqref{def:D-tail-ellip}.
\end{lemma}

\begin{proof}
     First, consider the weak formulation~\eqref{eq:weak_formulation} in $ [t_o - S, t_1]$ for some $t_1\in(t_o-S,t_o]$. 
By testing the weak formulation with $\boldsymbol{\tau}_{-h} \phi(x,t)$ and applying discrete integration by parts, we obtain
    \begin{align*}
        0=&\int_{\Omega} \phi\bth  u \,\dx \bigg|_{t_o-S}^{t_1}
        -
        \int_{t_o-S}^{t_1}\int_{\Omega} \partial_t\phi\bth  u \,\dx \dt \\
        &+ \int_{t_o-S}^{t_1}\iint_{\R^N\times\R^N} 
        \bigg( \Phi \bigg( \frac{u_{h}(x,t) - u_{h}(y,t)}{|x-y|^{s}} \bigg) - \Phi \bigg( \frac{u(x,t) - u(y,t)}{|x-y|^{s}} \bigg) \bigg) \nonumber \\
        &\phantom{==========}
        \cdot \frac{\phi(x,t)-\phi(y,t)}{\vert x-y \vert^{N+s}}   \,\dy \dx \dt, \nonumber
    \end{align*}
    where $u_h (x,t):= u(x+h,t)$. 
    We rewrite the diffusion term by introducing
    \begin{align*}
        K(x,y,t) =  \frac{\Phi \bigg( \dsty\frac{u_{h}(x,t) - u_{h}(y,t)}{|x-y|^{s}} \bigg) - \Phi \bigg( \frac{u(x,t) - u(y,t)}{|x-y|^{s}} \bigg)}{\dsty\frac{\bth u(x,t) - \bth u(y,t)}{| x-y |^{s}}} \frac{1}{|x-y|^{N+2s}}.
    \end{align*}
    Observe that $K : \mathbb{R}^{N} \times \mathbb{R}^{N} \times [t_o - S, t_o] \to \mathbb{R}$ is measurable and satisfies
    \begin{align}\label{ineq:kernel_estimate}
        \frac{\Lambda^{-1}}{ | x-y |^{N+2s}} \le K(x,y,t) = K(y,x,t) \le \frac{ \Lambda}{ | x-y |^{N+2s}},
    \end{align}
    due to~\eqref{assumption1:nonlinear_integrand_kernel}, for any $x,y \in \mathbb{R}^{N}$
    and $t \in [t_o - S, t_o]$. Thus, 
    \begin{align}
    \label{trump}
        0=&\int_{\Omega} \phi\bth  u \,\dx \bigg|_{t_o-S}^{t_1}
        -
        \int_{t_o-S}^{t_1}\int_{\Omega} \partial_t\phi\bth  u \,\dx \dt \\
        &+ 
        \int_{t_o-S}^{t_1}\iint_{\mathbb{R}^N\times\mathbb{R}^N} \big( \bth u(x,t) - \bth u(y,t) \big) \nonumber \\
        &\phantom{===========}
        \cdot \big( \phi(x,t)-\phi(y,t) \big) K(x,y,t) \,\dy \dx \dt. \nonumber
    \end{align}
In \eqref{trump} we choose $\phi(x,t) = \eta^2\bth  u(x,t)  $, where $\eta\in C^1(Q_{R,S}(z_o),[0,1])$
is a cut-off function satisfying 
\begin{equation*}
\left\{
    \begin{array}{c}
      \spt\eta\subset  Q_{\frac12 (R+r),S}(z_o),    \\
        \eta\equiv 1 \quad\text{in}\quad Q_{r,\tau} (z_o), \\
        \Vert D \eta \Vert_{L^{\infty}} \le \frac{4}{R-r}, \quad \vert \partial_t \eta \vert \le \frac{2}{S-\tau}.
    \end{array}\right.
\end{equation*} 
Standard calculations modulo a time mollification, cf.~\cite[Proof of Proposition~2.6]{DeFilippis} and \cite[Appendix~B]{Liao_24_2}, give
    \begin{align}\label{obama}
        \frac{1}{2} & \int_{B_R\times\{t_1\}} \eta^{2}\vert \bth u  \vert^{2}    \,\dx + \frac{1}{2} \int_{t_o-S}^{t_1} \int_{B_R}\partial_t \eta^{2} \vert \bth  u  \vert^2    \,\dx\dt \nonumber \\
        &\quad+\int_{t_o-S}^{t_1}\iint_{\mathbb{R}^N\times\mathbb{R}^N} \big( \bth u(x,t) - \bth u(y,t) \big) \nonumber \\
        &\phantom{===========}
        \cdot \big( \phi(x,t)-\phi(y,t) \big) K(x,y,t) \,\dy \dx \dt=0.
    \end{align}
    It suffices to deal with the third term on the left.
    Decomposing $\mathbb{R}^{N}$ into $B_R$ and its complement $\mathbb{R}^{N} \setminus B_R$ gives
    \begin{align*}
        \int_{t_o-S}^{t_1} & \int_{B_R} \int_{B_R} \big( \bth u(x,t) - \bth u(y,t) \big) \\
        &\phantom{=====}
        \cdot \big(  \bth u(x,t)  \eta^{2}(x,t) - \bth u(y,t)  \eta^{2}(y,t) \big) K(x,y,t) \,\dx \dy \dt \\
        &\phantom{=}
        + 2 \int_{t_o-S}^{t_1}\int_{B_{\frac12 (R+r)}}  \int_{\mathbb{R}^N \setminus B_R} \big( \bth u(x,t) - \bth u(y,t) \big) \\
        &\phantom{================}
        \cdot \bth u(x,t)  \eta^{2}(x,t) K(x,y,t) \,\dy \dx  \dt \\
        &=: \mathfrak{D} + 2\mathfrak{T}.
    \end{align*}
    We proceed by estimating the term $\mathfrak{D}$.    
    By applying~\eqref{ineq:kernel_estimate} to the kernel $K(x,y,t)$ and Lemma~\ref{lem:ineq_algebraic_a_b_c_d_e_f_p_euqal_to_2} for
    $a =  \bth u(x,t)$, $b=  \bth u(y,t)$, $e=\eta(x,t)$ and $f=\eta(y,t)$ to the integrand, we deduce 
    \begin{align*}
        \mathfrak{D} &\ge 
        \frac14
        \Lambda^{-1} \int_{t_o-S}^{t_1}\int_{B_R} \int_{B_R} \frac{ \vert \bth u(x,t) - \bth u(y,t)  \vert^2  ( \eta^{2}(x,t) + \eta^{2}(y,t)  ) }{\vert x-y \vert^{N+2s}} \,\dx\dy\dt \\
        &\phantom{=}
        - \frac12
        \Lambda^{-1} \int_{t_o-S}^{t_1}\int_{B_R} \int_{B_R} \frac{ ( \vert \bth u(x,t) \vert + \vert \bth u(y,t) \vert  )^{2}   }{\vert x-y \vert^{N+2s}}   \vert \eta(x,t) - \eta(y,t)  \vert^2 \,\dx\dy\dt \\
        &= : 
        \tfrac14
        \mathfrak{D}_1 -  
        \tfrac12
        \mathfrak{D}_2.
    \end{align*}
    Next, we further estimate $\mathfrak{D}_1$ from below by recalling that $\eta\equiv 1$ in $B_r\times (t_o-\tau,t_o]$:
    \begin{align*}
        \mathfrak{D}_1 &\ge \Lambda^{-1} \int_{t_o-\tau}^{t_o} \int_{B_r} \int_{B_r} \frac{\vert \bth u(x,t)  - \bth u(y,t)  \vert^{2} }{\vert x-y \vert^{N+2s}} \,\dx\dy\dt .
    \end{align*}
   For $\mathfrak{D}_{2}$, note that by the Lipschitz continuity of $\eta$, we obtain 
   \begin{align*}
      \vert \eta(x,t) - \eta(y,t) \vert^{2} \le  \Vert D \eta \Vert^{2}_{L^{\infty}} \vert x - y \vert^{2}  \le \frac{16}{(R-r)^2} \vert x - y \vert^2.
 \end{align*}
   Hence, due to symmetry, we deduce
   \begin{align*}
       \mathfrak{D}_2 
       &\le \frac{16 \Lambda^{-1}}{ (R-r)^2} \int_{t_o - S}^{t_o} \int_{B_R} \int_{B_R} \frac{ ( \vert \bth u(x,t) \vert + \vert \bth u(x,t) \vert  )^2}{\vert x - y \vert^{N+2s-2}} \,\dx\dy\dt \\
       &\le 
       \frac{64 \Lambda^{-1} }{ (R-r)^2} \int_{t_o-S}^{t_o} \int_{B_R}\bigg[\underbrace{\int_{B_R}
     \frac{ 1 }{|x-y|^{N-(1-s)2}}\dy}_{\le C(N,s) R^{2(1-s)}}\bigg] \vert \bth u(x,t) \vert^{2} \,\dx\dt \\
     &\le \frac{C}{  R^{2s}} \Big( \frac{R}{R - r} \Big)^{2} \int_{t_o - S}^{t_o} \int_{B_R} \vert \bth u \vert^{2} \,\dx\dt,
   \end{align*}
   where $C=C(\mathfrak{data})$.
   Finally, we treat the nonlocal term $\mathfrak{T}$. In particular, by applying~\eqref{ineq:kernel_estimate} to the kernel $K(x,y,t)$, Lemma~\ref{lem:first_tail_estimate} and Lemma~\ref{lem:prelim_tail_estimate} slice-wise it follows that 
   \begin{align*}
       \vert \mathfrak{T} \vert &\le \Lambda \int_{t_o-S}^{t_o}\int_{B_{\frac12 (R+r)}} \bigg\vert \int_{\mathbb{R}^N \setminus B_R} \frac{\bth u(x,t) - \bth u(y,t)}{\vert x-y\vert^{N+2s}} \,\dy \bigg\vert  \eta^{2}(x,t) |\bth u(x,t)|    \,\dx  \dt \\
       &\le \frac{C}{R^{2s}}\Big(\frac{R}{R-r}\Big)^{N+2s} \int_{t_o-S}^{t_o}\int_{B_{\frac12 (R+r)}}\eta^{2} \vert \bth  u \vert^{2}    \,\dx\dt \\
       &\phantom{=}
       + \frac{C \vert h \vert}{R^{2s+1}} 
    \Big(\frac{R}{R-r}\Big)^{N+2s+1} \int_{t_o-S}^{t_o} \mathcal T (t;B_{2R}) \int_{B_{\frac12 (R+r)}}  \eta^{2} |\bth u |  \, \dx \dt \\
    &=: \mathfrak{T}_{1} + \mathfrak{T}_{2},
    \end{align*}
    For the term $\mathfrak{T}_{1}$ we just recall that $0\le\eta\le1$. 
    We estimate $\mathfrak{T}_2$ by again taking advantage of the fact that $0\le\eta\le1$ and subsequently taking supremum in time of the spatial integral and then apply  Hölder's and Young's inequalities. To shorten the notation, we omit $B_{2R}$ from $\mathcal{T}$ and set $p=N+2s+1$. Then,
   \begin{align*}
       \mathfrak{T}_{2} 
       &\le  \frac{C \vert h \vert}{R^{2s+1}} \Big( \frac{R}{R-r} \Big)^p \bigg[ \int_{t_o - S}^{t_o} \mathcal{T}(t) \,\dt \bigg] \bigg[ \esssup_{t \in (t_o - S, t_o)} \int_{B_{R} \times \{t\}} \eta\vert \bth u  \vert \,  \dx \bigg] \nonumber \\
       &\le  \frac{C \vert h \vert}{R^{2s+1}} \Big( \frac{R}{R-r} \Big)^p \bigg[ \int_{t_o - S}^{t_o} \mathcal{T}(t) \,\dt \bigg] 
       |B_R|^{\frac12}\bigg[ \esssup_{t \in (t_o - S, t_o)} \int_{B_{R} \times \{t\}} \eta^{2}\vert \bth u  \vert^2 \,  \dx \bigg]^{\frac12} \nonumber \\
       &\le 
        \frac{C \vert h \vert^{2}}{R^{4s+2}} \Big( \frac{R}{R-r} \Big)^{2p} R^{N} \bigg[ \int_{t_o - S}^{t_o} \mathcal{T}(t) \,\dt \bigg]^{2}  
       + \frac14  \esssup_{t \in (t_o - S, t_o)}   \int_{B_R \times \{ t \} }\eta^{2} \vert \bth u \vert^2   \, \dx     . 
    \end{align*}
       Now, we combine all estimates in \eqref{obama}  and use the arbitrariness of $t_1 \in (t_o-\tau, t_o)$ to conclude the proof.
\end{proof}

The second kind of energy estimates concern the truncations of $\bth u$; it will be used to upgrade the obtained $W^{1,2}$-regularity to the Lipschitz regularity of weak solutions.
\begin{lemma}\label{lem:general_energy_inequality}
Let $s \in(0,1)$.  
If $u$ is a local weak solution to~\eqref{eq:PDE} in the sense of Definition~\ref{def:weak_solution}, then for any absolutely continuous function $k(\cdot)$ in $(0,T)$, any cylinder $Q_{R,S}(z_o)  \Subset \Omega_T$, any $h \in \mathbb{R}^{N} \setminus \{ 0 \}$ with $0 < \vert h \vert \le \frac34 R$ and any non-negative, Lipschitz function $\zeta(x,t)$ in $Q_{R,S}(z_o)$ such that $\spt\z(\cdot, t)\subset B_R(x_o)$ for any $t \in (t_o - S, t_o)$ 
we have
\begin{align*}
         \esssup_{t\in[t_o - S, t_o]}&\int_{B_R\times\{t\}}  \zeta^{2}w_{\pm}^{2}   \,\dx  + \Lambda^{-1}
         \iiint_{\mathbb{Q}_{R,S}} \frac{\vert  \z w_{\pm}(x,t) - \z w_{\pm}(y,t)  \vert^2}{\vert x-y \vert^{N+2s}}   \,\dy \dx  \dt \nonumber \\
          &\phantom{=}
         + \Lambda^{-1} \iiint_{\mathbb{Q}_{R,S}} \frac{\z^{2}(x,t) w_{\pm}(x,t) w_{\mp}(y,t)}{\vert x-y \vert^{N+2s}}   \,\dy \dx  \dt \nonumber \\
        &\le 2 \Lambda \iiint_{\mathbb{Q}_{R,S}} \max \{ w_{\pm}^{2}(x,t), w_{\pm}^{2}(y,t)  \}  \frac{ \vert \zeta(y,t) - \zeta(x,t) \vert^2}{\vert x-y \vert^{N+2s}} \,\dy \dx  \dt  \\
        &\phantom{=} +  4 \Lambda \iint_{Q_{R,S}}\big[\z^2 w_{\pm}^2(x,t) + |k(t)|\z^2 w_{\pm}(x,t) \big] \int_{\mathbb{R}^N \setminus B_{R}} \frac{\dy}{\vert x-y \vert^{N+2s}}  \, \dx  \dt
          \\
        &\phantom{=}+ 4 \Lambda \iint_{Q_{R,S}} \z^2w_{\pm}(x,t)\bigg|\int_{\mathbb{R}^N \setminus B_{R}} \frac{  \bth u(y,t)}{\vert x-y \vert^{N+2s}}  
         \,\dy  \bigg|\,\dx \dt   \\
        &\phantom{=} 
        \mp 2  \iint_{Q_{R,S}} k^{\prime} \zeta^{2}w_{\pm} \,\dx\dt 
        +  \iint_{Q_{R,S}} w_{\pm}^2 \vert \partial_t \zeta^{2} \vert \,\dx\dt \\
        &\phantom{=}
        + \int_{B_R\times\{t_o -S\}} \zeta^2 w_{\pm}^{2} \,\dx , \nonumber 
    \end{align*}
    where we defined $w:=\bth  u - k $ and omitted $x_o$ and $z_o$.  
\end{lemma}

\begin{proof}
    We will only deal with the case of $w_{+}$.
    Let $t_1\in(t_o-S,t_o]$ and take identity \eqref{trump} from the proof of Lemma~\ref{cor:second_energy_inequality} as starting point. 
    Next, we choose the test function 
    \begin{align*}
        \phi(x,t)
        :=
        \z^2(x,t)w_+(x,t)  
        \equiv 
        \z^2(\bth  u(x,t)  - k(t))_{+} ,
    \end{align*}
    which can be justified by a time mollification. Integrating by parts  
    the time part becomes
    \begin{align*}
        \int_{B_R} & \phi \bth  u  \,\dx \bigg\vert_{t_o - S}^{t_1}-\int_{t_o-S}^{t_1}\int_{B_R}\partial_t \phi \bth  u\,\dx \dt \nonumber \\
        &=
        \int_{t_o-S}^{t_1}\int_{B_R}\phi\partial_t(\bth  u) \,\dx \dt \nonumber \\
        &= 
        \frac{1}{2} \int_{B_R\times\{t_1\}}  \zeta^{2}w_{+}^{2} \,\dx - \frac{1}{2} \int_{B_R\times\{t_o - S\}}  \zeta^{2}w_{+}^{2} \,\dx  
 - 
        \frac{1}{2} \int_{t_o-S}^{t_1} \int_{B_R} \partial_t \zeta^{2} w_{+}^2  \,\dx\dt \nonumber \\
        &\phantom{=}
         +\int_{t_o-S}^{t_1}\int_{B_R} k^{\prime} \zeta^{2} w_{+}  \,\dx\dt. 
    \end{align*}
As for the diffusion part $\mathfrak{D}$, the last term in \eqref{trump}, we split the integral and take advantage of $\spt\zeta\subset B_{R} \times (t_o - S, t_o)$, as well as the integrand's symmetry, to get
    \begin{align*}
        \mathfrak{D}&= \int_{t_o-S}^{t_1}\int_{B_{R}}\int_{B_{R}} \big( \bth u(x,t) - \bth u(y,t) \big)  
         \big(\phi(x,t)-\phi(y,t)\big) K(x,y,t) \,\dy \dx  \dt \\
        &\phantom{=}
        + 2 \int_{t_o-S}^{t_1}\int_{B_{R}}\int_{\mathbb{R}^N \setminus B_{R}} \big( \bth u(x,t) - \bth u(y,t) \big)  
        \phi(x,t) K(x,y,t) \,\dy \dx  \dt \\
        &=: 
        \mathfrak{D}_{1} + \mathfrak{D}_{2}.
    \end{align*}
      For $\mathfrak{D}_1$, we use the identity
    \begin{align*}
        \big( \bth u(x,t) &- \bth u(y,t) \big) \big(  \phi (x,t) - \phi (y,t) \big) \nonumber\\
        &=  [ \zeta w_{+}(x,t) - \zeta w_{+}(y,t)   ]^2 + w_{-}(x,t) w_{+}(y,t) \zeta^2 (y,t) \\
        &\phantom{=}
        + w_{-}(y,t) w_{+}(x,t) \zeta^2 (x,t)
        - w_{+}(x,t) w_{+}(y,t)  [ \zeta(x,t) - \zeta(y,t)  ]^2. \nonumber
    \end{align*}
    Thus, using this identity 
    and applying~\eqref{ineq:kernel_estimate} to the kernel 
    we deduce  
   \begin{align*}
        \mathfrak{D}_{1} 
        &\ge 
        \int_{t_o-S}^{t_1}\int_{B_{R}}\int_{B_{R}} | \zeta  w_{+}(x,t) - \zeta  w_{+}(y,t) |^{2} K(x,y,t)    \,\dy \dx \dt \\
        &\phantom{=}
        + \int_{t_o-S}^{t_1}\int_{B_{R}}\int_{B_{R}} w_{-}(y,t) w_{+}(x,t) \zeta^2 (x,t) K(x,y,t)  \,\dy \dx  \dt \nonumber \\
        &\phantom{=}
        -  \int_{t_o-S}^{t_1}\int_{B_{R}}\int_{B_{R}}  w_{+}(x,t) w_{+}(y,t) | \zeta(y,t) - \zeta(x,t) |^{2} K(x,y,t)  \,\dy \dx \dt \nonumber \\
        &\ge \Lambda^{-1} \int_{t_o-S}^{t_1}\int_{B_{R}}\int_{B_{R}} \frac{| \zeta  w_{+}(x,t) - \zeta  w_{+}(y,t) |^{2}}{| x-y |^{N+2s}}   \,\dy \dx \dt \nonumber\\
        &\phantom{=}
        + \Lambda^{-1} \int_{t_o-S}^{t_1}\int_{B_{R}}\int_{B_{R}} \frac{w_{-}(y,t) w_{+}(x,t) \zeta^2 (x,t) }{| x-y |^{N+2s}}  \,\dy \dx  \dt \nonumber \\
        &\phantom{=}
        - \Lambda \int_{t_o-S}^{t_1}\int_{B_{R}}\int_{B_{R}} \max \{ w_{+}(x,t)^{2}, w_{+}(y,t)^{2}  \} \frac{ | \zeta(y,t) - \zeta(x,t) |^{2} }{| x-y |^{N+2s}}  \,\dy \dx \dt.
    \end{align*}
     Moreover, note that for the last term in the above display we used the estimate
    \begin{align*}
         w_{+}(x,t) w_{+}(y,t) \le \max \{ w_{+}^{2}(x,t), w_{+}^{2}(y,t)  \}.
    \end{align*} 
    As for $\mathfrak{D}_{2}$, we apply~\eqref{ineq:kernel_estimate} and split the term into
    \begin{align*}
        |\mathfrak{D}_{2}|
        &\le 2  \Lambda \bigg|\int_{t_o-S}^{t_1}\int_{B_{R}}\int_{\mathbb{R}^N \setminus B_{R}} \frac{\bth u(x,t)  }{\vert x-y \vert^{N+2s}}  
        \cdot\phi(x,t)  \,\dy \dx  \dt \bigg|\\
        &\quad +2  \Lambda \bigg|\int_{t_o-S}^{t_1}\int_{B_{R}}\int_{\mathbb{R}^N \setminus B_{R}} \frac{  \bth u(y,t)}{\vert x-y \vert^{N+2s}}  
        \cdot\phi(x,t)  \,\dy \dx \dt \bigg|\\
        &=:\mathfrak{D}_{2,1}+\mathfrak{D}_{2,2}.
    \end{align*}
    The first term is estimated simply by
    \begin{align*}
        \mathfrak{D}_{2,1} &\le 2  \Lambda \int_{t_o-S}^{t_1}\int_{B_{R}}\int_{\mathbb{R}^N \setminus B_{R}} \frac{\z^2 w_{+}^2(x,t) + |k(t)|\z^2 w_{+}(x,t) }{\vert x-y \vert^{N+2s}}  
          \,\dy \dx  \dt, 
    \end{align*}
    and the second term by
    \begin{align*}
        \mathfrak{D}_{2,2}\le 2  \Lambda \int_{t_o-S}^{t_1}\int_{B_{R}} \z^2 w_{+}(x,t)\bigg|\int_{\mathbb{R}^N \setminus B_{R}} \frac{  \bth u(y,t)}{\vert x-y \vert^{N+2s}}  
         \,\dy  \bigg|\,\dx \dt.
    \end{align*}
    Inserting these inequalities above and using the arbitrariness of $t_1\in(t_o-S,t_o]$ finishes the proof. 
\end{proof}

\section{Lipschitz regularity}

\subsection{$W^{1,2}$-regularity}
In this subsection we improve the preset $W^{s,2}$-regularity  of  local weak solutions to $W^{1,2}$. The approach consists of first exploiting the energy estimate for finite differences in Lemma~\ref{cor:second_energy_inequality} to obtain estimates for second order finite differences, and then performing a tailored finite-iteration argument to upgrade the spatial differentiability.

The statement we aim to prove is as follows.

\begin{proposition}\label{prop:W12}
Suppose that $s \in (0,1)$. 
    Then, for any locally bounded, local weak solution $u$ to~\eqref{eq:PDE} in the sense of Definition~\ref{def:weak_solution},
    we have that
    \begin{align*}
        u \in L_{\loc}^{2}(0,T; W_{\loc}^{1,2}(\Omega)).
    \end{align*}
    Moreover, for
    $Q_{2R}\equiv Q_{2R}(z_o) \Subset \Omega_T$,  
 there exists a constant $C=C(\mathfrak{data})$ such that the estimate
    \begin{align*}
        \biint_{Q_{R}} \vert D u \vert^{2} \,\dx\dt \le \frac{C}{R^{2}}\mathcal{T} ^{2}(Q_{2R}) 
    \end{align*}
    holds,
    where $\mathcal{T} (Q_{2R})$ is defined as in \eqref{def:D-tail-par}.
\end{proposition}
As mentioned above, before proceeding to the proof, we need a lemma concerning second-order finite differences.

\begin{lemma}\label{lem:second_order_finite_difference_bound_of_higher_order_derivatives_in_L_2}
    Suppose the assumptions of Proposition~\ref{prop:W12} hold true.
    Then,  for any $Q_{R}  \equiv Q_{R} (z_o) \Subset \Omega_T$   
    we have
    \begin{align*}
        \iint_{Q_{\frac14 R}} &\vert \bth^2  u \vert^{2} \,\dx\dt
        \le  C  \frac{|h|^{2s} }{R^{2s} } 
           \bigg[ \iint_{Q_R} \vert \bth  u \vert^{2} \,\dx\dt + \frac{\vert h \vert^{2}}{R^{2}} R^{N+2s}  \mathcal{T} ^{2}(Q_{2R})  \bigg]  
    \end{align*}
    for any $h \in \mathbb{R}^{N} \setminus \{ 0 \}$ with $\vert h \vert \le \frac18 R$, where $C=C(\mathfrak{data})$ and $\mathcal{T} (Q_{2R})$ is as in \eqref{def:D-tail-par}.
\end{lemma}

\begin{proof}
Upon rescaling, it suffices to prove it for $R=2$ and $z_o=(0,0)$.
 By applying Lemma~\ref{cor:second_energy_inequality}
 over $Q_1\subset Q_2$ 
    we deduce  
    that  
     \begin{align*}
          \int_{-1}^{0}[   \bth   u(\cdot,t)]_{W^{s,2}(B_{1})}^{2} \,\dt \le C  \iint_{Q_2} \vert \bth    u  \vert^2 \,\dx\dt 
         +C  \vert h \vert^{2}  \mathcal{T}^2 (Q_4)  .
     \end{align*}
Next, we apply Lemma~\ref{lem:ineq_for_finite_differences_of_sobolev_functions} slice-wise with $(v,\gamma,R)  = (  \bth   u ,  s, 1)$, integrate in time and use the above display to get
\begin{align*}
        \iint_{Q_{\frac12}} \vert \bth^2  u  \vert^{2} \,\dx \dt  &\le C \vert h \vert^{2s} \bigg[  \int_{-1}^{0}[   \bth   u(\cdot,t)]_{W^{s,2}(B_{1})}^{2} \,\dt + \iint_{Q_1} \vert \bth u \vert^{2} \,\dx\dt \bigg] \\
         &\le C  \vert h \vert^{2s}  
          \bigg[ \iint_{Q_{2}} \vert \bth  u \vert^{2} \,\dx\dt +\vert h \vert^{2} \mathcal{T}^2 (Q_4)  \bigg] \nonumber
    \end{align*}
    for $0 < \vert h \vert \le \frac12$, where $C=C(\mathfrak{data})$. 
\end{proof}

Now, we perform finite iteration and establish the aim of this subsection.
\smallskip

\noindent\textbf{Proof of Proposition~\ref{prop:W12}.}
By a proper scaling argument, it suffices to show the desired estimate with $R=1$.
    For $i\in \N_0$ set $\rho_{i} := 4^{-i}$ and $Q_i:=Q_{\rho_i}$; 
    we apply Lemma~\ref{lem:second_order_finite_difference_bound_of_higher_order_derivatives_in_L_2} with $R=\rho_{i} $  to get
\begin{align*} 
        \iint_{Q_{i+1}} \vert \bth^2  u \vert^{2} \,\dx\dt \le C    \vert h \vert^{2s}    \bigg[ \iint_{Q_i} \vert \bth    u  \vert^{2} \,\dx\dt + \vert h \vert^{2}  \mathcal{T} ^{2}(Q_{2\rho_i})   \bigg]  
    \end{align*}
    for some $C=C(i, \mathfrak{data})$ and for any $h \in \mathbb{R}^{N} \setminus \{ 0 \}$ with $\vert h \vert \le \frac18 \rho_{i}$, where
    \begin{align*} 
    \mathcal{T} ( Q_{2\rho_i}) := 
    \Vert  u \Vert_{L^{\infty}(Q_{2\rho_{i}})}
    +
    \mathrm{Tail}(u;   Q_{2\rho_{i}}).
    \end{align*}
    To simplify the above estimate, we note by Lemma~\ref{lem:prelim_tail_estimate} that
    \begin{align*}
      \mathcal{T} ( Q_{2\rho_i})  &\le C(i, N) \mathcal{T} ( Q_{2}). 
    \end{align*}
    Hence, there exists some $\widehat{C}_i=\widehat{C}_i(\mathfrak{data})$ such that
    \begin{align}\label{iteration-basis} 
        \iint_{Q_{i+1}} \vert \bth^2  u \vert^{2} \,\dx\dt &\le \widehat{C}_{i+1}  \vert h \vert^{2s}     \bigg[ \iint_{Q_i} \vert \bth    u  \vert^{2} \,\dx\dt +  \vert h \vert^{2} \mathcal{T} ^{2} ( Q_{2})  \bigg] \nonumber \\
        &\le  \widehat{C}_{i+1}  \vert h \vert^{2s} \mathcal{T} ^{2} ( Q_{2}) 
    \end{align}
    holds for all $i\in\N_0$ and for any $h \in \mathbb{R}^{N} \setminus \{ 0 \}$ with $\vert h \vert \le \frac18 \rho_{i+1}$. Next, we apply 
   Lemma~\ref{lem:domokos} 
    with $\gm=s$ joining \eqref{iteration-basis}$_{i=0}$,
    and conclude that
    \begin{align}\label{iteration-1}\nonumber
                \iint_{Q_{1}} \vert \bth  u \vert^{2} \,\dx\dt
                & \le \mathfrak{C}  \vert h \vert^{2s}     \bigg[\rho_o^{-2s}\iint_{Q_o} \vert   u  \vert^{2} \,\dx\dt+ \widehat{C}_1  \mathcal{T}^{2} ( Q_{2}) 
                \bigg] \nonumber \\  
                &\le  C_1 \vert h \vert^{2s}  \mathcal{T}^2 ( Q_{2}), 
    \end{align}
   where $\mathfrak{C}=\mathfrak{C}(s)$ and $C_1=\mathfrak{C}\widehat{C}_1\rho_o^{-2s}$, holds true for any $h \in \mathbb{R}^{N} \setminus \{ 0 \}$ with $\vert h \vert \le \frac{1}{8}\rho_{1}$.
    If $s\in[\frac12,1)$, we are done.
If $s\in[\frac14,\frac12)$, we first plug \eqref{iteration-1} back in \eqref{iteration-basis}$_{i=1}$ to get
    \begin{align}\label{iteration-2} 
        \iint_{Q_{2}} \vert \bth^2  u \vert^{2} \,\dx\dt &\le \widehat{C}_2  \vert h \vert^{2s}     \bigg[ C_1 \vert h \vert^{2s}  \mathcal{T} ^{2} ( Q_{2}) +  \vert h \vert^{2} \mathcal{T} ^{2} ( Q_{2})  \bigg]\\ \nonumber
        &\le 2\widehat{C}_2C_1 |h|^{2\cdot 2s} \mathcal{T}^{2} ( Q_{2})  ; 
    \end{align}
then, we apply Lemma~\ref{lem:domokos} with $\gm=2s<1$ joining \eqref{iteration-2} and conclude that
    \begin{align*}
                \iint_{Q_{2}} \vert \bth  u \vert^{2} \,\dx\dt& \le \mathfrak{C}  \vert h \vert^{2\cdot 2s}     \bigg[\rho_1^{-2s}\iint_{Q_1} \vert    u  \vert^{2} \,\dx\dt+ 2\widehat{C}_2C_1  \mathcal{T}^{2} ( Q_{2})  \bigg]\\
                &\le C_2 \vert h \vert^{2\cdot 2s}    \mathcal{T}^{2} ( Q_{2})
    \end{align*}
    for $C_2=2\mathfrak{C}\widehat{C}_2\rho_1^{-2s}C_1$ and for any $h \in \mathbb{R}^{N} \setminus \{ 0 \}$ with $\vert h \vert \le \frac{1}{8}\rho_{2}$.
Using induction,  if $s\in [\frac{1}{2i}, \frac{1}{2(i-1)} )$, $i\in\N$, then we will have
    \begin{align*}
                \iint_{Q_{i}} \vert \bth  u \vert^{2} \,\dx\dt
                \le C_i \vert h \vert^{i\cdot 2s}    \mathcal{T}^{2} ( Q_{2})
    \end{align*}
    for   $C_i=2\mathfrak{C}\widehat{C}_i\rho_{i-1}^{-2s}C_{i-1}$ and for any $h \in \mathbb{R}^{N} \setminus \{ 0 \}$ with $\vert h \vert \le \frac{1}{8}\rho_{i}$. Therefore, for $s\in(0,1)$, there must be $i_o\in\N$ such that $s\in [\frac{1}{2i_o}, \frac{1}{2(i_o-1)})$. The above estimate then gives that
    \begin{align*}
                \iint_{Q_{i_o}} \vert \bth  u \vert^{2} \,\dx\dt
                \le C_{i_o} \vert h \vert^{i_o\cdot 2s}   \mathcal{T} ^{2} ( Q_{2})\le C_{i_o} \vert h \vert^{2}    \mathcal{T} ^{2} ( Q_{2})
    \end{align*}
    for $C_{i_o}=(2\mathfrak{C})^{i_o}4^{si_o(i_o-1)} \prod_{j=1}^{i_o}\widehat{C}_j$ and for any $h \in \mathbb{R}^{N} \setminus \{ 0 \}$ with $\vert h \vert \le \frac{1}{8}\rho_{i_o}$. Note that $i_o$ is determined by $s$, and hence $C_{i_o}$ depends only on $N$, $s$ and $\Lambda$.
    Using Lemma~\ref{lem:bound_of_finite_difference_implies_bound_of_first_derivative} slice-wise, we can now conclude that 
        \begin{align*}
                \iint_{Q_{i_o}} \vert  D u \vert^{2} \,\dx\dt
               \le C     \mathcal{T}^{2} ( Q_{2})
    \end{align*}
    for some $C=C(\mathfrak{data})$. A covering argument will complete the proof.

\subsection{$W^{1,\infty}$-regularity}
With the $W^{1,2}$-regularity at hand, we aim to further upgrade the regularity to  $W^{1,\infty}$  in this subsection. The approach consists in performing De Giorgi-type iterations based on the energy estimates from Lemma~\ref{lem:general_energy_inequality}.

\begin{proposition}\label{prop:Lipschitz:1}
Let $s \in (0,1)$ 
and $u$
be a locally bounded, weak solution to~\eqref{eq:PDE} in the sense of Definition~\ref{def:weak_solution}.  
    Then,     
    \begin{align*}
        u \in L_{\loc}^{\infty}(0,T; W_{\loc}^{1,\infty}(\Omega)).
    \end{align*}
    Moreover, if $Q_{32R}\equiv Q_{32R}(z_o)   \Subset \Omega_T $, then  
    \begin{align*}
    \|D u \|_{L^\infty(Q_R)} 
         &\le \frac{C}{R}\Big[\|u\|_{L^\infty(Q_{8R} )} +  \mathrm{Tail}(u;Q_{8R} ) \Big] 
\end{align*}
holds true for some $C=C(\mathfrak{data}) > 1$.
\end{proposition}

\begin{proof} 
    We closely follow the lines of~\cite[Theorem 1.5]{Liao_24_1}. Suppose $u$ is a bounded, weak solution to~\eqref{eq:PDE} in $Q_{32}$.
 For $k>0 $ and $ n \in \mathbb{N}_{0}$  we introduce    
    \begin{align*}
	\left\{
	\begin{array}{c}
 k_n= k - 2^{-n} k, \quad  \rho_n=2^{-1} + 2^{-n-1},
\\[5pt]
 \widetilde{\rho}_n=\frac12(\rho_n+\rho_{n+1}), \quad \widehat{\rho}_n=\frac34\rho_n+ \frac14\rho_{n+1},\quad \overline{\rho}_n=\frac14\rho_n+ \frac34\rho_{n+1},\\[5pt]
  B_n=B_{\rho_n}, \quad \widetilde{B}_n=B_{\widetilde{\rho}_n},\quad \widehat{B}_n=B_{\widehat\rho_n},\quad \overline{B}_n=B_{\overline{\rho}_n},\\[5pt]
\dsty Q_n=Q_{\rho_n} ,\quad\widetilde{Q}_n=Q_{\widetilde{\rho}_n},\quad
\widehat{Q}_n=Q_{\widehat{\rho}_n},\quad\overline{Q}_n=Q_{\overline{\rho}_n}.
\end{array}
	\right.
\end{align*}
    Observe that  $Q_{n+1}\subset\overline{Q}_n\subset\widetilde{Q}_n\subset\widehat{Q}_n\subset Q_n$.
Choose a cutoff function $\z$ in $Q_n$, vanishing outside $\widehat{Q}_{n}$ and
equal to $1$ in $\widetilde{Q}_{n}$, such that
\begin{equation*}
|D\z|\le 2^{n+5} \quad\text{ and }\quad |\partial_t \z|\le  2^{n+5} .
\end{equation*}
    Next, we introduce 
    \begin{align*}
    \left\{
	\begin{array}{c}
        w(x,t) := (\bth u(x,t)-\mathsf{k}_{n+1}(t))_{+} , \quad \widetilde{w}(x,t) := ( \bth u(x,t)  - \mathsf{k}_n(t) )_{+},\\ [5pt]
        \displaystyle\mathsf{k}_{n}(t):=\ell(t) + k_{n},\quad 
        \ell(t):= \boldsymbol{\gm} \int_{ - 1}^{t}  \mathcal{T}(\tau; B_4)   \,\d\tau
        \end{array}\right.
    \end{align*}
    for  some  $\boldsymbol{\gm}>0$ to be specified, where we defined as in \eqref{def:D-tail-ellip}
    \begin{align*}
        \mathcal{T}(t;B_4) := \Vert u(t) \Vert_{L^{\infty}(B_4)} + \mathrm{Tail}(u(t); B_4).
    \end{align*}
    The energy inequality in Lemma~\ref{lem:general_energy_inequality} on $Q_n$ with the above choice of $\z$ and $w$ gives:
\begin{align}\label{newton} 
    \esssup\limits_{t \in [- \widetilde{\rho}^{2s}_{n},0]} &\int_{\widetilde{B}_n\times\{t\}} w^{2} \,\dx + \Lambda^{-1}\iiint_{\widetilde{\mathbb{Q}}_n} \frac{\vert w(x,t) - w(y,t) \vert^{2}}{\vert x-y \vert^{N+2s}} \,\dx\dy\dt \nonumber \\\nonumber
        &\le 2 \Lambda \iiint_{\mathbb{Q}_n}
        \max\{ w^{2}(x,t), w^{2}(y,t) \} \frac{  \vert \zeta(y,t) - \zeta(x,t)  \vert^2}{\vert x-y \vert^{N+2s}} \,\dy \dx  \dt  \\\nonumber
        &\phantom{=} +  4 \Lambda \iint_{Q_{n}}\big[\z^2 w^2(x,t) + \mathsf{k}_{n+1}(t)\z^2 w(x,t) \big] \int_{\mathbb{R}^N \setminus B_{n}} \frac{\dy}{\vert x-y \vert^{N+2s}}  \, \dx  \dt\\
        &\phantom{=} 
        + 4 \Lambda \iint_{Q_n} \z^2w(x,t)   \bigg\vert \int_{\mathbb{R}^N \setminus B_n} \frac{ \bth u(y,t)}{\vert x-y \vert^{N+2s}} \,\dy \bigg\vert \,\dx\dt \nonumber \\
        &\phantom{=} 
        - 2 \iint_{Q_n} \ell^{\prime} \zeta^{2} w \,\dx\dt 
        +  \iint_{Q_n}  w^2 \vert \partial_t \zeta^{2} \vert   \,\dx\dt \nonumber \\
&=:\mathbf{I}+\mathbf{II}+\mathbf{III}-\mathbf{IV}+\mathbf{V}. 
\end{align}
Here, the symbols $\widetilde{\mathbb{Q}}_n$  and $\mathbb{Q}_n$ should be self-explanatory. 
Estimating the first and fifth terms   is  straightforward, cf.~\cite{Liao_24_1}. Namely, 
  we have $\mathbf{I}$ and $\mathbf{V}$ bounded by
    \begin{align*}
          C 2^{2n} \iint_{Q_n} \widetilde{w}^{2}  \,\dx  \dt.
    \end{align*}
 For $\mathbf{II}$     we observe that $|x|\le \widehat{\rho}_{n}$ due to $\spt \z(\cdot,t)\subset \widehat{B}_n$ 
    and $\vert y \vert \ge \rho_{n}$. Then,
        \begin{align}\label{x-y}
        \frac{\vert y-x\vert}{\vert y \vert} \ge  1- \frac{\widehat{\rho}_{n}}{\rho_{n}} = \frac{1}{4} \bigg( \frac{\rho_{n} - \rho_{n+1}}{\rho_{n}} \bigg) \ge \frac{1   }{2^{n+4}} .
    \end{align} 
    By \eqref{x-y}, we have
    $$\mathbf{II}\le    C2^{n(N+2s)} \iint_{Q_{n}}\big[\z^2 w^2(x,t) + \mathsf{k}_{n+1}(t)\z^2 w(x,t) \big]  \, \dx  \dt.$$
    To continue, we recall $\widetilde{w}=(\bth u-\mathsf{k}_{n})_+$ and estimate
    \begin{align}\label{def:set_A_n}
        \vert A_n \vert &:= \vert \{ \bth u \ge \mathsf{k}_{n+1} \} \cap Q_n \vert \le \frac{2^{2(n+1)}}{k^2} \iint_{Q_n} \widetilde{w}^{2} \,\dx\dt ,
    \end{align}
    from which, applying Hölder's inequality and using $w  \le \widetilde{w} $, we deduce
    \begin{align}\label{w-less-w-tilde}
        \iint_{Q_n} w   \,\dx\dt \le \bigg( \iint_{Q_n} w^{2} \,\dx\dt \bigg)^{\frac{1}{2}} \vert A_n \vert^{\frac{1}{2}} \le \frac{2^{n+1}}{k} \iint_{Q_n}\widetilde{w}^{2}  \,\dx\dt.
    \end{align}
Therefore, upon using \eqref{w-less-w-tilde} and $\mathsf{k}_{n+1}(t)\le k+\ell(0)$, we have
    \begin{align*}
        \mathbf{II}  &\le   C 2^{(N+2s+1)n}   \Big[1+ \frac{\ell(0) }{k} \Big]\iint_{Q_n} \widetilde{w}^{2}  \,\dx\dt.
    \end{align*}
    Now, we will have a closer look at the third and forth terms, which need to be packed. 
    To begin with, 
    we decompose $\mathbf{III}$ as follows:
    \begin{align*}
        \mathbf{III} &\le 4 \Lambda \iint_{Q_n} \z^2w(x,t)   \bigg\vert \int_{\mathbb{R}^N \setminus B_2} \frac{ \bth u(y,t)}{\vert x-y \vert^{N+2s}} \,\dy \bigg\vert \,\dx\dt\\
        &\quad +4 \Lambda \iint_{Q_n} \z^2 w(x,t)   \bigg\vert \int_{B_2 \setminus B_n} \frac{ \bth u(y,t)}{\vert x-y \vert^{N+2s}} \,\dy \bigg\vert \,\dx\dt\\
        &=:\mathbf{III}_1+\mathbf{III}_2.
    \end{align*}
    For $\mathbf{III}_1$ we apply  Lemma~\ref{lem:first_tail_estimate} with $(r, R)=(1,2)$ and find that for any $0<|h|\le  \frac{1}{4}  $,
    \begin{align*}
        \mathbf{III}_1  
        \le   4 \Lambda C \vert h \vert  \iint_{Q_n} \mathcal{T}(t;B_4 )    \zeta^{2} w  \,\dx\dt\le   \Lambda C_*  \vert h \vert  \iint_{Q_n} \mathcal{T}(t;B_4 )    \zeta^{2} w  \,\dx\dt.
    \end{align*}
    Then, by definition of $\ell(t)$ we find that $\mathbf{III}_1 -\mathbf{IV}\le 0$, if we choose
    \begin{align*}
         2\boldsymbol{\gm} =  \Lambda C_*   \vert h \vert  .
    \end{align*}
    This choice of $\boldsymbol{\gm}$ also fixes
    \begin{align*}
         \ell(t) =   \tfrac12  \Lambda C_*   \vert h \vert   \int_{- 1}^{t} \mathcal{T}(\tau;  B_4 )  \,\d\tau. \nonumber
    \end{align*}
   For $\mathbf{III}_2$  we use \eqref{x-y} and \eqref{w-less-w-tilde} to estimate
    \begin{align*}
        \mathbf{III}_2 &\le  C 2^{(N+2s)n}  \iint_{Q_n} \z^2w(x,t)     \int_{B_2 \setminus B_n} \frac{ |\bth u(y,t)|}{\vert y \vert^{N+2s}} \,\dy   \,\dx\dt\\
        &\le  C 2^{(N+2s)n} \underbrace{\esssup_{t\in [-1,0]}\int_{B_{2}}|\bth u(y,t)|\,\dy}_{=:\mathfrak{I}} \iint_{Q_n} \z^2w        \,\dx\dt\\
        &\le C 2^{(N+2s+1)n} \frac{\mathfrak{I}}{k}\iint_{Q_n} \widetilde{w}^2        \,\dx\dt.
    \end{align*}
    Therefore, we have
    \begin{align*}
        \mathbf{III}-\mathbf{IV}\le \mathbf{III}_{1}+\mathbf{III}_{2} -\mathbf{IV}\le\mathbf{III}_{2} \le C 2^{(N+2s+1)n} \frac{\mathfrak{I}}{k}\iint_{Q_n} \widetilde{w}^2        \,\dx\dt,
    \end{align*}
    and, up to now we have obtained
    \begin{align*}
    \mathbf{I}+\mathbf{II}+\mathbf{III}-\mathbf{IV}+\mathbf{V} \le
          C 2^{(N+2s+1)n}  \bigg[1 + \frac{\ell(0)}{k} + \frac{\mathfrak{I}}{k}  \bigg] \iint_{Q_n} \widetilde{w}^2   \,\dx\dt.
    \end{align*}
   Let us stipulate that
   \begin{align}\label{ineq:k_requirement_1_in_Lipschitz_proof}
       \ell(0) +\mathfrak{I}  \le k
   \end{align}
   such that the terms in the bracket are bounded by $2$.
    Collecting all the estimates on the right-hand side of~\eqref{newton} gives us 
    \begin{align}\label{ineq:short_energy_estimate}
         \esssup\limits_{t \in [ - \widetilde{\rho}^{2s}_{n},0]}& \int_{\widetilde{B}_n\times\{t\}} w^{2} \,\dx + \Lambda^{-1} \iiint_{\widetilde{\mathbb{Q}}_n} \frac{\vert w(x,t) - w(y,t) \vert^{2}}{\vert x-y \vert^{N+2s}} \,\dx\dy\dt \nonumber \\
        &\le C 2^{(N+2s+1)n}  \iint_{Q_n} \widetilde{w}^2   \,\dx\dt.
    \end{align}
    At this point consider a cutoff function $\phi\in C^1(\widetilde{Q}_n,[0,1])$, which vanishes outside of $\overline{Q}_n $, is equal to $1$ in $Q_{n+1} $ and satisfies $\vert D \phi \vert \le  2^{n+4} $. By applying  Hölder's inequality and the embedding from Proposition~\ref{prop:parabolic_embedding}, we obtain
    \begin{align}\label{DG-Sobolev}
        \iint_{Q_{n+1}}& w^{2} \,\dx\dt \le \iint_{\widetilde{Q}_n}  ( \phi  w )^2 \,\dx\dt \le \bigg[ \iint_{\widetilde{Q}_n} (  \phi w )^{\frac{2(N+2s)}{N}} \,\dx\dt \bigg]^{\frac{N}{N+2s}} \vert A_n \vert^{\frac{2s}{N+2s}} \nonumber\\
        &\le C \bigg[  \iiint_{\widetilde{\mathbb{Q}}_n} \frac{\vert  \phi w(x,t)   -  \phi w(y,t)   \vert^{2}}{\vert x - y \vert^{N+2s}} \,\dx\dy\dt
        +
        \iint_{\widetilde{Q}_n}  (  w \phi  )^2 \,\dx\dt \bigg]^{\frac{N}{N+2s}} \nonumber \\
        &\phantom{==}
        \cdot \bigg[ \esssup\limits_{t \in [ - \widetilde{\rho}^{2s}_n, 0]} \fint_{\widetilde{B}_n}  (  w \phi  )^2 \,\dx \bigg]^{\frac{2s}{N+2s}} \vert A_n \vert^{\frac{2s}{N+2s}}, 
    \end{align}
    for $A_n$ being defined according to~\eqref{def:set_A_n}.
    The numerator of the integrand inside the fractional norm is estimated by
    \begin{align*}
        \vert \phi w(x,t) &- \phi w(y,t) \vert^{2} \nonumber \\
        &\le 2 \phi^{2}(x,t)\vert  w(x,t) -  w(y,t) \vert^{2}  + 2 w^{2}(y,t) \vert \phi(x,t) - \phi(y,t) \vert^{2}.
    \end{align*}
    Therefore, we deduce that
    \begin{align*}
          &\iiint_{\widetilde{\mathbb{Q}}_n}  \frac{\vert \phi  w(x,t) - \phi  w(y,t) \vert^{2}}{\vert x - y \vert^{N+2s}} \,\dx\dy\dt \nonumber \\
        &\le 2  \iiint_{\widetilde{\mathbb{Q}}_n} \frac{\vert  w(x,t) -  w(y,t) \vert^{2}}{\vert x - y \vert^{N+2s}} \,\dx\dy\dt 
        + 2   \iiint_{\widetilde{\mathbb{Q}}_n} \frac{w^{2}(y,t)\vert \phi(x,t) - \phi(y,t) \vert^{2}}{\vert x - y \vert^{N+2s}} \,\dx\dy\dt \nonumber\\
        &\le 2   \iiint_{\widetilde{\mathbb{Q}}_n} \frac{\vert  w(x,t) -  w(y,t) \vert^{2}}{\vert x - y \vert^{N+2s}} \,\dx\dy\dt 
        +C 2^{2n} \iint_{\widetilde{Q}_n} w^{2} \,\dy\dt \nonumber\\
        &\le C 2^{(N+2s+1)n}  \iint_{Q_n} \widetilde{w}^2   \,\dx\dt. \nonumber
    \end{align*}
    Note that we used~\eqref{ineq:short_energy_estimate} for the last estimate, where the constant $\Lambda^{-1}$ was brought to the right-hand side and absorbed into $C$, as well as $w \le \widetilde{w}$ and $\widetilde{Q}_{n}  \subset Q_n $. Similarly, the sup-in-time term in \eqref{DG-Sobolev} can be bounded by using~\eqref{ineq:short_energy_estimate}, whereas for $\vert A_n \vert$ we use~\eqref{def:set_A_n}.
    In this way, we deduce from \eqref{DG-Sobolev} that
\begin{align}\label{Einstein}
        \iint_{Q_{n+1}}  w^{2} \,\dx\dt  
        &\le C\bigg[  2^{(N+2s+1)n}      \iint_{Q_n} \widetilde{w}^{2}  \,\dx  \dt  
        + 
        \iint_{Q_n} \widetilde{w}^{2} \,\dx\dt \bigg]^{\frac{N}{N+2s}} \nonumber\\
        &\phantom{==}
        \cdot \bigg[   2^{(N+2s+2)n}     \iint_{Q_n} \widetilde{w}^{2}  \,\dx  \dt  \bigg]^{\frac{2s}{N+2s}} 
        \cdot \bigg[\frac{2^{2(n+1)}}{k^2} \iint_{Q_n} \widetilde{w}^{2} \,\dx\dt \bigg]^{\frac{2s}{N+2s}} \nonumber\\
        &\le  \frac{C b^{n} }{  k^{2 \cdot \frac{2s}{N+2s}}} \bigg[ \iint_{Q_n} \widetilde{w}^{2} \,\dx\dt \bigg]^{1+\frac{2s}{N+2s}}
\end{align}
for some $b=b(N)>1$.
    Defining
    \begin{align*}
        Y_{n}:= \frac{1}{k^2}\iint_{Q_n}  ( \bth u  - \ell(t) - k_n  )_{+}^{2} \,\dx\dt
    \end{align*}
     and recalling the definitions of $w$ and $\widetilde{w}$ will enable us as rewrite~\eqref{Einstein} as
    \begin{align*}
        Y_{n+1} \le 
        C b^{n}    Y_{n}^{1+\frac{2s}{N+2s}}  .
    \end{align*}
     This encodes a fast geometric convergence, cf.~\cite[Chapter I, Lemma 4.1]{DiBenedetto}. In fact, if
    \begin{align}\label{ineq:requirement_2}
        Y_{o}&= 
        \frac{1}{k^{2}}  
        \iint_{Q_o} (  \bth u  - \ell(t) )_{+}^{2} \,\dx\dt \le C^{-\frac{N+2s}{2s}} b^{-(\frac{N+2s}{2s})^2}
    \end{align}
      then
    \begin{align*}
        \lim_{n \to \infty} \iint_{Q_n}(  \bth u  - \ell(t) - k_n )_{+}^{2} \,\dx\dt = 0, 
    \end{align*}
    which amounts to
    \begin{align*}
         \bth u  \le \ell + k \quad \text{a.e. in }    Q_\frac12.
    \end{align*}
    We properly choose $k$ to fulfill both  \eqref{ineq:k_requirement_1_in_Lipschitz_proof}
    and~\eqref{ineq:requirement_2}, and obtain from the last line that
    \begin{align}\label{ineq:intermediate_local_boundedness_estimate}
         \esssup\limits_{  Q_\frac12}   \bth u  &\le \mathfrak{I} +  C   \vert h \vert \int_{ - 1}^{0} \mathcal{T}(t; B_4 )  \,\d t  + C\bigg[ \iint_{Q_1} |\bth u|^{2} \,\dx\dt \bigg]^{\frac{1}{2}} \nonumber \\
        &=: \mathfrak{I} + \mathfrak{L}_{1} + \mathfrak{L}_{2}.
    \end{align}

    Our next aim is to estimate the right-hand side in the above display  in terms of $|h|$.
    First, we will estimate $\mathfrak{I}$ by  the energy inequality from Lemma~\ref{cor:second_energy_inequality}. Indeed, we omit the second term on the left of the energy inequality and choose proper cut-off functions to get       
    \begin{align*}
        \mathfrak{I}&=    \esssup\limits_{t\in [  -1, 0 ]}  \int_{B_{2}} \vert \bth u \vert  \,\dx 
        \le C \bigg[ \esssup\limits_{t\in [ -1, 0]}  \int_{B_{2}} \vert \bth u \vert^{2}  \,\dx \bigg]^{\frac{1}{2}}\\
        &\le C
        \bigg[\iint_{Q_{4}} \vert \bth u\vert^{2} \,\dx\dt\bigg]^{\frac{1}{2}} + 
          C \vert h \vert  \bigg[ \int_{-2}^{0} \mathcal{T}(t;B_4) \,\dt \bigg]  \\
        &=:  \mathfrak{I}_{1} + \mathfrak{I}_{2}  ,
    \end{align*}
    where we defined as in \eqref{def:D-tail-ellip}
    \begin{align*}
        \mathcal{T}(t;B_4) := \Vert u(t) \Vert_{L^{\infty}(B_{4})} + \mathrm{Tail}(u(t); B_4).
    \end{align*} 
In order to estimate $\mathfrak{I}_{1}$, we apply Proposition~\ref{prop:W12} to get  that
\begin{align*}
    \mathfrak{I}_{1} \le  C \vert h \vert  \mathcal{T} ,
\end{align*}
where we defined as in \eqref{def:D-tail-par}
\begin{align*}
    \mathcal{T}\equiv \mathcal{T}(Q_{8}) &:=   \Vert u  \Vert_{L^{\infty}(Q_{8})} +  \mathrm{Tail}(u;Q_8).
\end{align*}
By Lemma~\ref{lem:prelim_tail_estimate} it is straightforward to see that
    $\mathcal{T} (t;B_4) \le C \mathcal{T} (t;B_8)$.
Consequently,
\begin{align*}
    \mathfrak{I}_{2} \le   C \vert h \vert     \mathcal{T} ,
\end{align*}
and, therefore, we get that
\begin{align*}
    \mathfrak{I} \le \mathfrak{I}_{1}+\mathfrak{I}_{2}\le   C \vert h \vert  \mathcal{T}.
\end{align*}
    The other two terms in \eqref{ineq:intermediate_local_boundedness_estimate} will be bounded by the same quantity. Indeed,
    by Lemma~\ref{lem:prelim_tail_estimate} it is not hard to see that
    \begin{align*}
        \mathfrak{L}_{1}\le  C\vert h \vert  \int_{ -2}^{0} \mathcal{T}(t;B_8) \,\dt  \le   C\vert h \vert   \mathcal{T}.
    \end{align*}
    Moreover, by applying Proposition~\ref{prop:W12}   we obtain that 
\begin{align*}
    \mathfrak{L}_{2}&\le C  \vert h \vert  \mathcal{T}.
\end{align*}
    Combining all the estimates from above yields the desired inequality
\begin{align}\label{ineq:esssup_bound_for_finite_differnce_in_Lipschitz_proof}
    \esssup\limits_{  Q_{\frac12}}  \bth u  
         &\le C |h| \mathcal{T},
\end{align}
for any $h \in \R^{N} \setminus \{ 0 \}$, as long as $|h|\le\frac14$.
To complete the proof, let us fix $(x,t)\in Q_\frac14$ and consider $y$ satisfying $|x-y|\le\frac14$. Apply the above estimate \eqref{ineq:esssup_bound_for_finite_differnce_in_Lipschitz_proof} with $h = y-x$ to deduce
\begin{align*}
      u( y,t ) - u(x,t)   
         &\le C |y-x| \mathcal{T}. 
\end{align*}
Switching the roles of $x$ and $y$ we may deduce that
\begin{align}\label{ineq:abs_value_of_finite_difference_quotient_bound}
     | u( y,t ) - u(x,t) | 
         &\le C |y-x| \mathcal{T}
\end{align}
holds for any $(x,t)\in Q_\frac14$ and $y$ satisfying $|x-y|\le\frac14$.
In other words, we have obtained
\begin{align}\label{ineq:final_estimate_in_Lipschitz_proof}
        \|D u \|_{L^\infty(Q_\frac14)} 
         &\le C \Big[\|u\|_{L^\infty(Q_{8} )} +  \mathrm{Tail}(u;Q_{8} ) \Big].
\end{align}
    A proper scaling argument finishes the proof.
\end{proof}

\subsection{Proof of Lipschitz estimates}
    First, 
    by Proposition~\ref{prop:Lipschitz:1} for $R=1$,~\cite[Theorem 1.6]{Liao_24_1} with $b=0$ and the fact that $u \in  L^1_\loc \big(0,T;L^1_{2s}(\R^N)\big)$ we obtain the estimate
    \begin{align*}
        \| D u \|_{L^\infty(Q_1)} &\le  C \Big[\|u\|_{L^\infty(Q_{8} )} +  \mathrm{Tail}(u;Q_{8} ) \Big] \\
         &\le C \bigg[ \int_{-32^{2s}}^{0} \int_{\mathbb{R}^{N}} \frac{u_{+}(x,t)}{1+|x|^{N+2s}} \,\dx\dt + \mathrm{Tail}(u;Q_{8} ) \bigg] \nonumber\\
         &\le C  \int_{-32^{2s}}^{0} \int_{\mathbb{R}^{N}} \frac{| u(x,t) |}{1+|x|^{N+2s}} \,\dx\dt. \nonumber
    \end{align*}
    We have proven the Lipschitz estimate in Theorem~\ref{Thm:main} after proper scaling and covering.

\section{Expansion of positivity}
Fix $h\in \R^N\setminus\{0\}$ and consider
\begin{align}\label{def:abbrevaition_for_diff_quotient}
    v \equiv v_h:=\frac{\bth u}{|h|},
\end{align}
defined on $\Om_T$. When no ambiguity arises, we omit $h$ from $v_h$.
Next, we state a reformulation of Lemma~\ref{lem:general_energy_inequality}. In fact, dividing both sides of the energy inequality in Lemma~\ref{lem:general_energy_inequality} by $|h|^2$, relabeling $k/|h|$ by $k$ and taking the supremum over the term involving $f$, defined as in \eqref{def:widetilde_f_as_nonlocal_term}, leads to the following energy inequality.

\begin{lemma}\label{lem:localized_energy_inequality} 
Let $s \in(0,1)$.  
If $u$ is a local weak solution to~\eqref{eq:PDE} in the sense of Definition~\ref{def:weak_solution}, then for any 
absolutely continuous function $k(\cdot)$ in $(0,T)$, any cylinder $Q_{R,S}(z_o) $ in $ \Omega_T$, any $h \in \mathbb{R}^{N} \setminus \{ 0 \}$ with $0 < \vert h \vert \le \frac34 R$ and any non-negative, Lipschitz function $\zeta(x,t)$ in $Q_{R,S}(z_o)$ such that $\spt\z(\cdot, t)\subset B_R(x_o)$ for any $t \in (t_o - S, t_o)$, 
we have  
\begin{align*}
    \esssup_{t\in[t_o - S, t_o]}&\int_{B_R \times\{t\}}  \zeta^{2} \widetilde{w}_{\pm}^{2}   \,\dx  + \Lambda^{-1}
         \iiint_{\mathbb{Q}_{R,S}} \frac{\vert  \z \widetilde{w}_{\pm}(x,t) - \z \widetilde{w}_{\pm}(y,t)  \vert^2}{\vert x-y \vert^{N+2s}}   \,\dy \dx  \dt \nonumber \\
         &\phantom{=}
         + \Lambda^{-1} \iiint_{\mathbb{Q}_{R,S}} \frac{\z^{2}(x,t) \widetilde{w}_{\pm}(x,t) \widetilde{w}_{\mp}(y,t)}{\vert x-y \vert^{N+2s}}   \,\dy \dx  \dt \nonumber \\
        &\le 2 \Lambda \iiint_{\mathbb{Q}_{R,S}} \max \{ \widetilde{w}_{\pm}^{2}(x,t), \widetilde{w}_{\pm}^{2}(y,t)  \}  \frac{ \vert \zeta(y,t) - \zeta(x,t) \vert^2}{\vert x-y \vert^{N+2s}} \,\dy \dx  \dt  \\
        &\phantom{=} +  4 \Lambda \iint_{Q_{R,S}}\big[\z^2 \widetilde{w}_{\pm}^2(x,t) + | k(t) | \z^2 \widetilde{w}_{\pm}(x,t) \big] \int_{\mathbb{R}^N \setminus B_{R}} \frac{\dy}{\vert x-y \vert^{N+2s}}  \, \dx  \dt
          \\
        &\phantom{=}+ 4 \Lambda \int_{t_o - S}^{t_o} \Vert f(\cdot, t) \Vert_{L^{\infty}( \spt\z(\cdot, t))} \int_{B_R} \z^2 \widetilde{w}_{\pm}(x,t) \,\dx \dt   \\
        &\phantom{=} 
        \mp 2  \iint_{Q_{R,S}} k^{\prime} \zeta^{2} \widetilde{w}_{\pm} \,\dx\dt 
        +  \iint_{Q_{R,S}} \widetilde{w}_{\pm}^2 \vert \partial_t \zeta^{2} \vert \,\dx\dt \\
        &\phantom{=}
        + \int_{B_R\times\{t_o-S\}} \zeta^2 \widetilde{w}_{\pm}^{2} \,\dx , \nonumber 
\end{align*}
where we omitted $x_o$, $z_o$ and defined $\widetilde{w}:= v - k $ as well as  
\begin{align}\label{def:widetilde_f_as_nonlocal_term}
    f(x,t):=  \int_{\mathbb{R}^N \setminus B_{R}} \frac{  v(y,t)}{\vert x-y \vert^{N+2s}}  
         \,\dy \in L^1_{\loc}(0,T; L^{\infty}_{\loc}(B_R)) .
\end{align}
\end{lemma}

With the energy estimates at hand, we aim to derive a modulus of continuity of $v\equiv v_h$ that is uniform in $h$. The program is a tailored version of \cite{Liao_24_2, Liao_24_3}. The gist is that proper energy estimates encode all information needed to perform regularity estimates.

 Introduce a reference cylinder $\mathcal{Q}:=B_{R} (x_o) \times (t_1, t_2] \subset \Omega_T$.
Assume that $\boldsymbol{\mu}^{\pm}$ and $\boldsymbol{\om}$ satisfy 
\begin{align}\label{ineq:condition_on_m_pm}
    \boldsymbol{\mu}^{+} := \esssup_{\mathcal{Q}} v, \quad \boldsymbol{\mu}^{-} := \essinf_{\mathcal{Q}} v, \quad \boldsymbol{\om} \ge  \boldsymbol{\mu}^{+} - \boldsymbol{\mu}^{-}.
\end{align}
Our goal in this section is to use the energy estimate from Lemma~\ref{lem:localized_energy_inequality} to prove the following expansion of positivity for $v$.

\begin{proposition}[Expansion of positivity]\label{prop:expansion_of_positivity}
    Let $u$ be a locally bounded, local weak solution to~\eqref{eq:PDE} in $\Omega_T$
    and suppose that for some $t_o \in (t_1, t_2)$ and some $\lambda,\, \xi \in (0,1)$ 
    we have
    \begin{equation*}
    \begin{array}{cc}
              |  \{ \pm  ( \boldsymbol{\mu}^{\pm} - v(\cdot, t_o)  ) \ge \xi \boldsymbol{\om}  \} \cap B_{\rho}(x_o)  | \ge \lambda | B_{\rho}|.  
    \end{array}
\end{equation*}
    Then, there exist constants $\theta,\, \eta \in (0,1)$, $\widetilde{\boldsymbol{\gamma}}$ depending only on the $\mathfrak{data}$ and $\lambda$, such that
    either
    \begin{align*}
    \widetilde{\boldsymbol{\gamma}} \int_{t_1}^{t_2} \Vert f(\cdot,t) \Vert_{L^{\infty}(B_{8 \rho }(x_o))} \, \dt>  \eta \xi \boldsymbol{\om}
\end{align*}
or
    \begin{align*}
        \pm( \boldsymbol{\mu}^{\pm} - v ) \ge \eta \xi \boldsymbol{\om} \quad \text{a.e. in}\>\> B_{2 \rho}(x_o) \times  ( t_o + \tfrac{1}{2}\theta \rho^{2s}, t_o + \theta \rho^{2s} ],
    \end{align*}
    provided 
    \begin{align*}
        B_{8 \rho} (x_o) \times ( t_o, t_o + \theta \rho^{2s} ] \subset \mathcal{Q}.
    \end{align*}
    Moreover, we have the dependencies $\theta \approx \lambda^{N+3}$ and $\eta \approx \lambda^{q}$ for some $q = q(N) >1$. 
\end{proposition}

\subsection{Preparatory tools}
Three lemmata will be prepared before we prove the expansion of positivity.
The first one is a De Giorgi-type lemma. Recall $\mathcal{Q}$, $\boldsymbol{\mu}^{\pm}$ and $\boldsymbol{\om}$ satisfying~\eqref{ineq:condition_on_m_pm}. Here and in what follows, we will use the cylinder $Q_{\rho}(\theta):=B_\rho(x_o)\times (t_o-\theta\rho^{2s}, t_o]$ and omit the vertex $(x_o,t_o)$.

\begin{lemma}[De Giorgi-type lemma]\label{lem:a_de_giorgi_type_lemma}
Let $u$ be a locally bounded, local weak solution to \eqref{eq:PDE} in $\Omega_T$.
    For $\theta,\, \xi \in (0,1)$, $\rho \in (0, \frac18 R)$  and $Q_{\rho}(\theta) 
    \subset \mathcal{Q}$, there exists a constant $\widetilde{\boldsymbol\gamma} >1$ depending only on the $\mathfrak{data}$ and $\nu \in (0,1)$ depending only on the $\mathfrak{data}$ and $\theta$, such that if 
\begin{equation*}
    \begin{array}{cc}
              |  \{ \pm  ( \boldsymbol{\mu}^{\pm} - v  ) \le \xi \boldsymbol{\om}  \} \cap Q_{\rho}(\theta)  | \le \nu | Q_{\rho}(\theta) |, 
    \end{array}
\end{equation*}
then either
\begin{align*}
    \widetilde{\boldsymbol{\gamma}} \int_{t_1}^{t_2} \Vert f(\cdot,t) \Vert_{L^{\infty}(B_{2 \rho })} \,\dt >  \xi \boldsymbol{\om}
\end{align*}
or
\begin{align*}
    \pm ( \boldsymbol{\mu}^{\pm} - v ) \ge \tfrac{1}{4} \xi \boldsymbol{\om} \quad \text{a.e. in } Q_{\frac{1}{2} \rho}(\theta).
\end{align*}
Moreover, we have the dependence $\nu \approx \theta^q$ for some $q=q(N) >1$.
\end{lemma}

\begin{proof}
     We prove the statement for $w_{-}$ with $\boldsymbol{\mu}^{-} = 0$. 
Upon a translation, we may assume $z_o = (0,0)$ and still use $\mathcal{Q}:=B_{R} \times (t_1, t_2] $ to denote the shifted reference cylinder.
Next, we introduce
    \begin{align*}
	\left\{
	\begin{array}{c}
 k_n = (\frac{1}{2} + \frac{1}{2^{n+1}})\xi \boldsymbol{\om}, \quad  \rho_n= (\frac{1}{2} + \frac{1}{2^{n+1}})\rho,
\\[5pt]
 \widetilde{\rho}_n=\frac12(\rho_n+\rho_{n+1}), \quad \widehat{\rho}_n=\frac34\rho_n+ \frac14\rho_{n+1},\quad \overline{\rho}_n=\frac14\rho_n+ \frac34\rho_{n+1},\\[5pt]
  B_n=B_{\rho_n}, \quad \widetilde{B}_n=B_{\widetilde{\rho}_n},\quad \widehat{B}_n=B_{\widehat\rho_n},\quad \overline{B}_n=B_{\overline{\rho}_n},\\[5pt]
\dsty Q_n=B_n \times ( - \theta \rho_{n}^{2s} ,0),\quad\widetilde{Q}_n=\widetilde{B}_n \times ( - \theta \widetilde{\rho}_{n}^{2s},0),\\[5pt]
\widehat{Q}_n=\widehat{B}_n \times ( - \theta \widehat{\rho}_{n}^{2s},0),\quad\overline{Q}_n=\overline{B}_n \times ( - \theta \overline{\rho}_{n}^{2s},0).
\end{array}
	\right.
\end{align*}
Note that  $Q_{n+1}\subset\overline{Q}_n\subset\widetilde{Q}_n\subset\widehat{Q}_n\subset Q_n$.
Let $\z$ in $Q_n$ be a cutoff function, vanishing outside $\widehat{Q}_{n}$ and
equal to $1$ in $\widetilde{Q}_{n}$ such that
\begin{align*}
    \vert D \zeta \vert \le \frac{2^{n+4}}{\rho} \quad \text{and} \quad \vert \partial_t \zeta \vert \le \frac{2^{2s(n+4)}}{\theta \rho^{2s}}.
\end{align*}
Applying the energy estimate from Lemma~\ref{lem:localized_energy_inequality} in $Q_n$ with $\zeta$ and 
\begin{align}\label{def:definition_of_w_minus-and_ell}
    \left\{
	\begin{array}{c}
        \widetilde{w}(x,t) :=  ( v(x,t)-\mathsf{k}_{n}(t)  ) , 
        \\ [5pt]
        \displaystyle\mathsf{k}_{n}(t):= k_{n} - \ell(t),\quad 
        \ell(t):= \widetilde{\boldsymbol\gamma} \int_{t_1}^{t} 
         \Vert f(\cdot, \tau) \Vert_{L^{\infty}(B_{2 \rho })} \,
        \d\tau, 
        \end{array}\right.
    \end{align}
 for $\widetilde{\boldsymbol\gamma} >0$ to be specified later,
    gives us
    \begin{align*}
        &\esssup_{t\in[- \theta \widetilde{\rho}_{n}^{2s}, 0]} \int_{\widetilde{B}_n \times\{t\}} \zeta^2 \widetilde{w}_{-}^2 \,\dx + \Lambda^{-1} \iiint_{\widetilde{\mathbb{Q}}_n}
        \frac{ \vert \zeta \widetilde{w}_{-}(x,t) - \zeta \widetilde{w}_{-}(y,t) \vert^{2}}{\vert x-y \vert^{N+2s}}  \,\dy\dx\dt \nonumber \\
        &\le 2 \Lambda \iiint_{\mathbb{Q}_n} \max \{ \widetilde{w}_{-}^{2}(x,t), \widetilde{w}_{-}^{2}(y,t)  \}  \frac{ \vert \zeta(y,t) - \zeta(x,t) \vert^2}{\vert x-y \vert^{N+2s}} \,\dy \dx  \dt  \\
        &\phantom{=} +  4 \Lambda \iint_{Q_{n}}\big[\z^2 \widetilde{w}_{-}^2(x,t) + | \displaystyle\mathsf{k}_{n+1}(t) | \z^2 \widetilde{w}_{-}(x,t) \big] \int_{\mathbb{R}^N \setminus B_{n}} \frac{\dy}{\vert x-y \vert^{N+2s}}  \, \dx  \dt
          \\
        &\phantom{=}+ 4 \Lambda \int_{- \theta \rho_{n}^{2s}}^{0} \Vert f(\cdot, t) \Vert_{L^{\infty}(\widehat{B}_n)} \int_{B_n} \z^2 \widetilde{w}_{-}(x,t) \,\dx \dt   \\
        &\phantom{=} 
        - 2  \iint_{Q_{n}} \ell^{\prime} \zeta^{2} \widetilde{w}_{-} \,\dx\dt 
        +  \iint_{Q_{n}} \widetilde{w}_{-}^2 \vert \partial_t \zeta^{2} \vert \,\dx\dt \nonumber   \\
        & =:\mathbf{I}+\mathbf{II}+\mathbf{III}-\mathbf{IV}+\mathbf{V}.
\end{align*}
Note that the last term from Lemma~\ref{lem:localized_energy_inequality} vanishes due to the support of $\zeta$.
Estimating the terms $\mathbf{I}$ and $\mathbf{V}$ is standard, that is,
\begin{align*}
    \mathbf{I}  \le C 2^{2n}  \frac{(\xi \boldsymbol{\om})^2}{\rho^{2s}} | A_n |
\quad\text{and}\quad
    \mathbf{V} \le \frac{2^{2sn}}{\theta \rho^{2s}} ( \xi \boldsymbol{\om} )^2 | A_n |,
\end{align*}
where
$$A_n := \{ v(x,t) \le  k_{n} - \ell(t)   \} \cap Q_n.$$
Next, we treat $\mathbf{II}$ as follows. First, observe that if $| y | \ge \rho_n$ and $| x | \le \widehat{\rho}_n$, then
\begin{align*}
    \frac{| y-x |}{|y|} \ge 1 - \frac{\widehat{\rho}_n}{\rho_n} = \frac{1}{4} \bigg( \frac{\rho_n - \rho_{n+1}}{\rho_n} \bigg) \ge \frac{1}{2^{n+4}}.
\end{align*}
Therefore, we obtain assuming $\ell(t_2)\le \xi\boldsymbol{\om}$ that
\begin{align*}
    \mathbf{II} &\le \frac{C 2^{(N+2s)n}}{\rho^{2s}} \iint_{Q_n} \big[\z^2 \widetilde{w}_{-}^2(x,t) + | \displaystyle\mathsf{k}_{n+1}(t) | \z^2 \widetilde{w}_{-}(x,t) \big] \,\dx\dt \\
    &\le \frac{C 2^{(N+2s)n}}{\rho^{2s}} ( \xi \boldsymbol{\om})^2 | A_n |,
\end{align*}
since $| \displaystyle\mathsf{k}_{n+1}(t) | \le | \displaystyle\mathsf{k}_{n+1} | + | \ell(t) | \le 2 \xi \boldsymbol{\om}$. For the third term, we have 
\begin{align*}
    \mathbf{III} &\le 4 \Lambda \int_{- \theta \rho_{n}^{2s}}^{0} \Vert f(\cdot, t) \Vert_{L^{\infty}(B_{2 \rho})} \int_{B_n} \z^2 \widetilde{w}_{-}(x,t) \,\dx \dt
\end{align*}
and then choose $\widetilde{\boldsymbol\gamma} = 2 \Lambda$ in the definition of $\ell(t)$, cf.~\eqref{def:definition_of_w_minus-and_ell}, such that $$\ell^{\prime}(t) = 2 \Lambda \Vert f(\cdot,t) \Vert_{L^{\infty}(B_{2 \rho })},$$ and hence $\mathbf{III} - \mathbf{IV} \le 0$. To summarize, we have deduced that
\begin{align*}
    &\esssup_{t\in[- \theta \widetilde{\rho}_{n}^{2s}, 0]} \int_{\widetilde{B}_n \times\{t\}} \zeta^2 \widetilde{w}_{-}^2 \,\dx + \Lambda^{-1} \iiint_{\widetilde{\mathbb{Q}}_n}
        \frac{ \vert \zeta \widetilde{w}_{-}(x,t) - \zeta \widetilde{w}_{-}(y,t) \vert^{2}}{\vert x-y \vert^{N+2s}}  \,\dy\dx\dt \nonumber \\
        &\le \frac{C 2^{(N+2s)n}}{\theta \rho^{2s}} ( \xi \boldsymbol{\om})^2 | A_n |.
\end{align*}
Proceeding  as in~\eqref{DG-Sobolev} -- \eqref{Einstein} we obtain  by defining $Y_n := \frac{| A_n |}{|Q_n|}$ that
\begin{align*}
    Y_{n+1} \le C \theta^{-  \frac{N}{N+2s}} b^{n} 
    Y_{n}^{1+\frac{2s}{N+2s}}.
\end{align*}
Hence, by the fast geometric convergence, cf.~\cite[Chapter I, Lemma 4.1]{DiBenedetto}, there exists a $\nu$ depending only on the $\mathfrak{data}$ and $\theta$, such that $Y_n \to 0$, if we require $Y_o \le \nu$, i.e.
\begin{align*}
     | \{ v(x,t) + \ell(t) \le  \xi \boldsymbol{\om} \} \cap Q_{\rho}(\theta)  | \le \nu | Q_{\rho}(\theta) |.
\end{align*} 
Then, $Y_n \to 0$ implies that
\begin{align*}
    v(x,t) + \ell(t) \ge   \tfrac{1}{2} \xi \boldsymbol{\om}  \quad \text{a.e. in}\>\> Q_{\frac{1}{2} \rho}(\theta).
\end{align*}
This implies that, 
if $\ell(t_2) \le \frac{1}{4} \xi \boldsymbol\omega$, by adjusting $\widetilde{\boldsymbol\gamma}$ appropriately, and
\begin{align*}
     |  \{ v(x,t)  <  \xi \boldsymbol{\om} \} \cap Q_{\rho}(\theta)  | \le \nu | Q_{\rho}(\theta) |,
\end{align*}
then
\begin{align*}
    v \ge  \tfrac{1}{4} \xi \boldsymbol{\om}  \quad \text{a.e. in } Q_{\frac{1}{2} \rho}(\theta).
\end{align*}
We conclude the proof by redefining $4 \widetilde{\boldsymbol\gamma}$ to $\widetilde{\boldsymbol\gamma}$.
\end{proof}

Next, we formulate a lemma, which propagates measure theoretical information forward in time.

\begin{lemma}\label{lem:propagate_measure_theoretical_information_forward_in_time}
Let $u$ be a locally bounded, local weak solution to~\eqref{eq:PDE} in $\Omega_T$.
Introduce the parameters $\xi, \lambda \in (0,1)$. There exist $\theta, \epsilon \in (0,1)$ depending only on the $\mathfrak{data}$ and $\lambda$, such that if for some $t_o \in (t_1, t_2)$ we have
\begin{equation*}
    \begin{array}{cc}
              |  \{ \pm  ( \boldsymbol{\mu}^{\pm} - v(\cdot, t_o)  ) \ge \xi \boldsymbol{\om}  \} \cap B_{\rho}(x_o)  | \ge \lambda | B_{\rho}|,  
    \end{array}
\end{equation*}
then either
\begin{align*}
    \frac{1}{\theta} \int_{t_1}^{t_2} \Vert f(\cdot,t) \Vert_{L^{\infty}(B_{2 \rho })} \,\dt >  \xi \boldsymbol{\om}
\end{align*}
or
    \begin{align*}
         |  \{ \pm  ( \boldsymbol{\mu}^{\pm} - v(\cdot,t)   ) \ge \epsilon \xi \boldsymbol{\om}  \} \cap B_{\rho}(x_o)  | \ge \tfrac{1}{2}\lambda | B_{\rho} | \quad \text{for all}\>\> t \in  ( t_o, t_o + \theta \rho^{2s}  ],
    \end{align*}
    provided $B_{\rho}(x_o) \times (t_o, t_o + \theta \rho^{2s} ] \subset \mathcal{Q}$. Moreover, we have $\epsilon \approx \lambda$ and $\theta \approx \lambda^{N+3}$.
\end{lemma}

\begin{proof}
    We will assume that $\boldsymbol{\mu}^{-}=0$, take $(x_o, t_o) = (0,0)$, and still denote the reference cylinder $B_R\times(t_1,t_2]$. Suppose that 
    \begin{align}\label{ineq:rewritten_first_alternative_forward_time_propagation}
        \frac{1}{\theta} \int_{t_1}^{t_2} \Vert f(\cdot,t) \Vert_{L^{\infty}(B_{2 \rho })} \,\dt  
        \le \xi \boldsymbol{\om}.
    \end{align}
    Let $\widetilde{w}_{-}(x,t) := ( v(x,t)  - \xi \boldsymbol{\om})_{-}$, consider $\alpha \in (0,1)$, which will be determined later, and choose a cutoff function $\zeta(x,t) \equiv \zeta(x)$, such that  $\z=1$ on $B_{(1-\alpha)\rho}$ and $\z=0$ outside  $B_{(1- \frac{\alpha}{2})\rho}$ satisfying $| D \zeta | \le \frac{4}{\alpha \rho}$. Hence,   the energy inequality from Lemma~\ref{lem:localized_energy_inequality} in $B_{\rho} \times (0, \theta \rho^{2s}]$ with $\widetilde{w}_{-}$,  $\zeta$ and $k = \xi \boldsymbol{\om}$ gives for all $t \in (0, \theta \rho^{2s})$ that
    \begin{align*}
        &\int_{B_{(1-\alpha) \rho} \times \{ t \}}  \widetilde{w}_{-}^{2} \,\dx \\
        &\le 2 \Lambda \int_{0}^{\theta \rho^{2s}} \int_{B_{\rho}} \int_{B_{\rho}} \max \{ \widetilde{w}_{-}^{2}(x,t), \widetilde{w}_{-}^{2}(y,t)  \}  \frac{ \vert \zeta(y,t) - \zeta(x,t) \vert^2}{\vert x-y \vert^{N+2s}} \,\dy \dx  \dt  \\
        &\phantom{=} +  4 \Lambda \int_{0}^{\theta \rho^{2s}} \int_{B_{\rho}} \big[\z^2 \widetilde{w}_{-}^2(x,t) + \xi \boldsymbol{\om} \z^2 \widetilde{w}_{-}(x,t) \big] \int_{\mathbb{R}^N \setminus B_{\rho}} \frac{\dy}{\vert x-y \vert^{N+2s}}  \, \dx  \dt
          \\
        &\phantom{=}+  4 \Lambda \int_{0}^{\theta \rho^{2s}} \Vert f(\cdot, t) \Vert_{L^{\infty}(\spt \z)} \int_{B_{\rho}} \z^2 \widetilde{w}_{-}(x,t) \,\dx \dt +  \int_{B_\rho \times \{ 0 \} } \widetilde{w}_{-}^{2}\,\dx  \\
        &=: \mathbf{I}+\mathbf{II}+\mathbf{III} + \mathbf{IV}.
    \end{align*}
    In the above application of Lemma~\ref{lem:localized_energy_inequality} we omitted the second and third terms on the left-hand side, whereas the last two terms on the right-hand side vanished because $k$ and $\zeta$ are independent of time. We estimate the term $\mathbf{I}$ by 
    \begin{align*}
        \mathbf{I} &\le 2 \Lambda \theta \rho^{2s} (\xi \boldsymbol{\om})^2 \int_{B_{\rho}} \int_{B_{\rho}} \frac{| \zeta(x) - \zeta(y) |^2}{| x-y |^{N+2s}} \,\dy\dx \\
        &\le 32 \Lambda \theta \rho^{2s} \frac{(\xi \boldsymbol{\om})^2}{(\alpha \rho)^2} \int_{B_{\rho}} \int_{B_{\rho}} \frac{1}{| x-y |^{N+2(s-1)}} \,\dy\dx \\
        &\le C \frac{\theta ( \xi \boldsymbol{\om})^2}{\alpha^2} | B_{\rho} |.
    \end{align*}
    For the second term, note that since $x \in \spt \zeta \subset B_{(1-\frac{\alpha}{2})\rho}$ and $y \in \mathbb{R}^{N} \setminus B_{\rho}$, we have
    \begin{align*}
        \frac{| x-y |}{|y|} \ge 1 - \frac{|x|}{|y|} \ge \frac{\alpha}{2}.
    \end{align*}
    Thus, it follows that
    \begin{align*}
        \mathbf{II} &\le \frac{C}{\alpha^{N+2s} \rho^{2s}} \int_{0}^{\theta \rho^{2s}} \int_{B_{\rho}} \big[\z^2 \widetilde{w}_{+}^2(x,t) + \xi \boldsymbol{\om} \z^2 \widetilde{w}_{+}(x,t) \big]   \, \dx  \dt \\
        &\le \frac{C}{\alpha^{N+2s} \rho^{2s}} | Q_{\rho}(\theta) | (\xi \boldsymbol{\om})^2 = C \frac{\theta (\xi \boldsymbol{\om})^2 }{\alpha^{N+2s} } | B_{\rho}| .
    \end{align*}
    Next, due to~\eqref{ineq:rewritten_first_alternative_forward_time_propagation} we obtain for $\mathbf{III}$ the estimate  
    \begin{align*}
        \mathbf{III} &\le C \xi \boldsymbol{\om} | B_{ \rho} | \int_{0}^{\theta \rho^{2s}} \Vert f(\cdot, t) \Vert_{L^{\infty}(B_{ 2 \rho })} \, \dt \\
        &\le C \theta (\xi \boldsymbol{\om})^2 | B_{\rho} |.
    \end{align*}
    Furthermore, we obtain for the term $\mathbf{IV}$ that
    \begin{align*}
        \mathbf{IV} \le (1- \lambda) (\xi \boldsymbol\omega)^2 | B_\rho |.
    \end{align*}
     By combining all the estimates from above we obtain for all $t \in (0, \theta \rho^{2s})$ that
    \begin{align}\label{ineq:very_short_energy_estimate_for_slicewise_De_giorgi}
        \int_{B_{(1-\alpha) \rho} \times \{ t \}}  \widetilde{w}_{-}^{2}  \,\dx \le (\xi \boldsymbol\omega)^2 \bigg[ (1-\lambda) + C \frac{\theta }{\alpha^{N+2s}} \bigg] | B_{\rho} |.
    \end{align}
    Let us estimate  the left-hand side of the previous display from below by
    \begin{align*}
        \int_{B_{(1-\alpha) \rho} \times \{ t \}}  \widetilde{w}_{-}^{2} \,\dx &\ge \int_{B_{(1-\alpha) \rho} \times \{ t \}}  \widetilde{w}_{-}^{2} \chi_{ \{ v <  \epsilon \xi \boldsymbol{\om} \} } \,\dx \\
        &\ge (1- \epsilon )^2 ( \xi \boldsymbol{\om})^2 | A_{\epsilon, (1-\alpha) \rho} (t) |,
    \end{align*}
    where
    \begin{align*}
        A_{\epsilon, (1-\alpha) \rho} (t) := \{ v(\cdot,t) <  \epsilon \xi \boldsymbol{\om} \} \cap B_{(1-\alpha) \rho}. 
    \end{align*}
    Now, observe that
    \begin{align}\label{ineq:estimate_of_level_sets}
        | A_{\epsilon, \rho} (t) | &= | A_{\epsilon, (1-\alpha) \rho} (t) \cup ( A_{\epsilon,  \rho} (t)  \setminus  A_{\epsilon, (1-\alpha) \rho} (t) ) | \nonumber \\
        &\le | A_{\epsilon, (1-\alpha) \rho} (t) | + | B_{\rho} \setminus B_{(1-\alpha)\rho} | \nonumber \\
        &\le  | A_{\epsilon, (1-\alpha) \rho} (t) | + N \alpha | B_{\rho} |.
    \end{align}
    Combining~\eqref{ineq:very_short_energy_estimate_for_slicewise_De_giorgi} and~\eqref{ineq:estimate_of_level_sets} yields for all $t \in (0, \theta \rho^{2s})$ that
    \begin{align*}
        | A_{\epsilon, \rho} (t) | &\le \frac{1}{(1-\epsilon)^2} \bigg[  (1-\lambda) + C \frac{\theta}{\alpha^{N+2}}  \bigg] | B_{\rho} | + N \alpha  | B_{\rho} |.
    \end{align*}
    In order to finish the proof we  choose the parameters $\alpha, \theta$ and $\epsilon$ in such a way that
    \begin{align*}
        N \alpha = \tfrac{1}{4} \lambda, \quad C \frac{\theta}{\alpha^{N+2}} = \tfrac{1}{8} \lambda \quad \text{and} \quad 
        \frac{1 - \frac{7}{8} \lambda}{(1- \epsilon)^2} \le 1 - \tfrac34 \lambda
    \end{align*}
    are satisfied. Therefore, we deduce for all $t \in (0, \theta \rho^{2s})$ that
    \begin{align*}
        | A_{\epsilon, \rho} (t) | \le  (1- \tfrac{1}{2} \lambda  ) | B_{\rho} |.
    \end{align*}
    This is just  rephrasing the desired estimate.
\end{proof}

We conclude this section by deriving a measure shrinking lemma.

\begin{lemma}\label{lem:measure_shrinking_lemma}
Let $u$ be a locally bounded, local weak solution to~\eqref{eq:PDE} in $\Omega_T$,
as well as $\sigma\in(0,\frac12)$ and $\theta, \xi, \lambda\in (0,1)$. Suppose that for all $t \in (t_o - \theta \rho^{2s}, t_o] \subset(t_1,t_2)$  
we have
    \begin{equation*}
    \begin{array}{cc}
             | \{ \pm ( \boldsymbol{\mu}^{\pm} - v(\cdot, t) ) \ge \xi \boldsymbol{\om} \} \cap B_{\rho}(x_o) | \ge \lambda | B_{\rho}|.  
    \end{array}
\end{equation*}
There exists $C>1$ depending only on the $\mathfrak{data}$, such that
either
\begin{align*}
    \frac{1}{\theta} \int_{t_1}^{t_2} \Vert f(\cdot,t) \Vert_{L^{\infty}(B_{2 \rho })} \,\dt >  \sigma \xi \boldsymbol{\om}
\end{align*}
or
\begin{align*}
    | \{ \pm ( \boldsymbol{\mu}^{\pm} - v ) \le \sigma \xi \boldsymbol{\om} \} \cap Q_{\rho}(\theta) | \le C \frac{\sigma}{\theta \lambda} | Q_{\rho}(\theta) |
\end{align*}
provided $Q_{4 \rho} (\theta) \subset \mathcal{Q}$.
\end{lemma}

\begin{proof}
    We take $(x_o, t_o) =(0,0)$, still denote the reference cylinder $B_R\times(t_1,t_2]$, and  prove the assertion for $\boldsymbol{\mu}^{-} = 0$.
    Next, introduce the truncation
    \begin{align*}
        \widetilde{w}_{-}(x,t) := ( v(x,t)  - \sigma \xi \boldsymbol{\om} )_{-}
    \end{align*}
    and a  cutoff function $\zeta$ in $B_{2 \rho}$, such that $\z=1$ in $B_{\rho}$ and $\z=0$ outside $B_{2 \rho}$ satisfying $| D \zeta | \le \frac{4}{\rho}$. Then, from Lemma~\ref{lem:localized_energy_inequality}, applied in $B_{4 \rho} \times (- \theta \rho^{2s}, 0]$ with $\widetilde{w}_{-}, \zeta$ and $k = \sigma \xi \boldsymbol{\om}$, we deduce
    \begin{align*}
        \Lambda^{-1} &\iint_{Q_{\rho}(\theta)} \widetilde{w}_{-}(x,t) \bigg( \int_{B_{4 \rho}} \frac{  \widetilde{w}_{+}(y,t)}{|x-y|^{N+2s}} \,\dy \bigg) \dx\dt \\
        &\le  2 \Lambda \int_{- \theta \rho^{2s}}^{0} \int_{B_{4 \rho}} \int_{B_{4 \rho}} \max\{ \widetilde{w}_{-}^{2}(x,t) , \widetilde{w}_{-}^{2}(y,t) \}   \frac{ \vert \zeta(x,t) - \zeta(y,t) \vert^{2}}{\vert x-y \vert^{N+2s}} \,\dy\dx\dt  \\
       &\phantom{=}
       +  4 \Lambda \int_{-\theta \rho^{2s}}^{0} \int_{B_{4\rho}} \big[\z^2 \widetilde{w}_{-}^2(x,t) + \sigma \xi \boldsymbol{\om} \z^2 \widetilde{w}_{-}(x,t) \big] \int_{\mathbb{R}^N \setminus B_{4\rho}} \frac{\dy}{\vert x-y \vert^{N+2s}}  \, \dx  \dt
          \\
        &\phantom{=}+  4 \Lambda \int_{-\theta \rho^{2s}}^{0} \Vert f(\cdot, t) \Vert_{L^{\infty}(\spt \zeta)} \int_{B_{4 \rho}} \z^2 \widetilde{w}_{-}(x,t) \,\dx \dt  \\
        &\phantom{=}
        + \int_{B_{4 \rho}\times\{x,  - \theta \rho^{2s}\}} \zeta^2 \widetilde{w}_{-}^{2} \,\dx \\
        &=: \mathbf{I} + \mathbf{II} + \mathbf{III} + \mathbf{IV}.
    \end{align*}
    Since, the terms $\mathbf{I}, \mathbf{II}$ and $\mathbf{III}$ are treated similar as in the proof of Lemma~\ref{lem:propagate_measure_theoretical_information_forward_in_time}, we will only present the main difference. Starting with the first term we obtain
    \begin{align*}
        \mathbf{I} 
        &\le C \theta \rho^{2s} \frac{ (\sigma \xi \boldsymbol{\om})^2}{\rho^2} \frac{| B_{4 \rho} |}{\rho^{2(s-1)}} \le C \frac{(\sigma \xi \boldsymbol{\om} )^2}{\rho^{2s}} | Q_{\rho}(\theta) |.
    \end{align*}
    Proceeding with $\mathbf{II}$ note that for $|x| \le \frac{3}{4} \rho$ and $|y| \ge 2 \rho$ we have
    \begin{align*}
        \frac{|x-y|}{|y|} \ge 1 - \frac{|x|}{|y|} \ge \frac{1}{4}.
    \end{align*}
    Thus, we obtain analogously as in the proof of Lemma~\ref{lem:a_de_giorgi_type_lemma} that
    \begin{align*}
        \mathbf{II} 
        &\le C \theta ( \sigma \xi \boldsymbol{\om} )^2 | B_{4 \rho} | \le C 4^{N} \frac{( \sigma \xi \boldsymbol{\om} )^2}{\rho^{2s}} | Q_{\rho}(\theta) |.
    \end{align*}
    For $\mathbf{III}$, since 
    $$\frac{1}{\theta} \int_{t_1}^{t_2} \Vert f(\cdot,t) \Vert_{L^{\infty}(B_{2 \rho })} \,\dt \le  \sigma \xi \boldsymbol{\om},$$ it follows that
    \begin{align*}
        \mathbf{III} 
        &\le C \theta (\sigma \xi \boldsymbol{\om})^2  | B_{4 \rho} | \le C 4^{N} \frac{( \sigma \xi \boldsymbol{\om} )^2}{\rho^{2s}} | Q_{\rho}(\theta) |.
    \end{align*}
    Next, $\mathbf{IV}$ is estimated by
    \begin{align*}
        \mathbf{IV} \le ( \sigma \xi \boldsymbol\omega)^2 | B_{4 \rho} | \le C \frac{( \sigma \xi \boldsymbol\omega)^2}{\theta \rho^{2s}} | Q_{\rho}(\theta) |.
    \end{align*}
    Note that we used $| B_{4 \rho} | = 4^{N} | B_{\rho} |$.
    By combining all these estimate, while shifting the constant $\Lambda^{-1}$ to the right-hand side and absorbing it into $C$, we can find the bound
    \begin{align}\label{ineq:upper_bound_for_left_middle_energy_term}
        \iint_{Q_{\rho}(\theta)} &\int_{B_{4 \rho}} \frac{ \widetilde{w}_{-}(x,t) \widetilde{w}_{+}(y,t)}{|x-y|^{N+2s}} \,\dy\dx\dt  \le C \frac{( \sigma \xi \boldsymbol{\om} )^2}{\rho^{2s}} | Q_{\rho}(\theta) |.
    \end{align}
    Note that $C$ only depends on the $\mathfrak{data}$. Next, by using the measure theoretical information the left-hand side can be estimated from below by
    \begin{align}\label{ineq:lower_bound_for_left_middle_energy_term}
        \iint_{Q_{\rho}(\theta)} &\widetilde{w}_{-}(x,t) \chi_{\{ v(x,t) \le  \frac{1}{2} \sigma \xi \boldsymbol{\om} \}} \bigg( \int_{B_{4 \rho}} \frac{   \widetilde{w}_{+}(y,t) \chi_{\{ v(y,t) \ge  \sigma \xi \boldsymbol{\om} \}}}{|x-y|^{N+2s}} \,\dy \bigg) \,\dx\dt \nonumber \\
        &\ge \tfrac{1}{2} \sigma \xi \boldsymbol{\om} | \{ v \le  \tfrac{1}{2} \sigma \xi \boldsymbol{\om} \} \cap Q_{\rho}(\theta) | \bigg( \frac{\xi \boldsymbol{\om} \lambda |B_{\rho}|}{2 ( 8 \rho)^{N+2s}} \bigg) \nonumber \\
        &= 2^{-(4N +2+4s)} \frac{(\xi \boldsymbol{\om})^2 \lambda \sigma}{\rho^{2s}} | \{ v \le  \tfrac{1}{2} \sigma \xi \boldsymbol{\om} \} \cap Q_{\rho}(\theta) |.
    \end{align}
We combine~\eqref{ineq:upper_bound_for_left_middle_energy_term} with~\eqref{ineq:lower_bound_for_left_middle_energy_term} and adjust the constants appropriately to finish the proof.
\end{proof}

\subsection{Proof of Proposition~\ref{prop:expansion_of_positivity}}

 We prove the assertion for $\boldsymbol{\mu}^{-} =0$ and $(x_o, t_o) = (0,0)$. The reference cylinder $B_R\times(t_1, t_2]$ will be used. Let us rewrite the measure theoretical information at the initial time $t_o =0$ in the larger ball $B_{4 \rho}$,  replace $\lambda$ by $4^{-N} \lambda$, and
apply Lemma~\ref{lem:propagate_measure_theoretical_information_forward_in_time} in order to obtain $\theta,\, \epsilon \in (0,1)$ depending only on the $\mathfrak{data}$ and $\lambda$, such that
\begin{align*}
    | \{  v(\cdot, t)  \ge  \epsilon \xi \boldsymbol{\om} \} \cap B_{4 \rho} | \ge \tfrac{1}{2}4^{-N}\lambda  | B_{4 \rho}| \quad \text{for all } t \in (0, \theta (4 \rho)^{2s}],
\end{align*}
provided we enforce
\begin{align*}
    \frac{1}{\theta} \int_{t_1}^{t_2} \Vert f(\cdot,t) \Vert_{L^{\infty}(B_{8 \rho })} \,\dt \le  \xi \boldsymbol{\om} .
\end{align*}
This slice-wise measure theoretical estimate enables us to apply Lemma~\ref{lem:measure_shrinking_lemma} in the cylinders 
\begin{align*}
    (0, \overline{t}) + Q_{4 \rho} ( \tfrac{1}{2} \theta ) := B_{4 \rho} \times  ( \overline{t} - \tfrac{1}{2}\theta (4 \rho)^{2s}, \overline{t}  ]
\end{align*}
for any arbitrary $\overline{t} \in  ( \frac{1}{2} \theta( 4 \rho)^{2s}, \theta (4 \rho)^{2s}  ]$ and with $\xi$, $\lambda$ replaced by $\epsilon \xi$, $\frac{1}{2} 4^{-N} \lambda$ respectively. 
Hence, the backward cylinders considered are contained in $B_{4 \rho} \times  ( 0, \theta (4 \rho)^{2s}  ]$. 
Now let $\nu$ be determined according to Lemma~\ref{lem:a_de_giorgi_type_lemma} depending only on the $\mathfrak{data}$ and $\theta$. Moreover, choose $\sigma$ according to Lemma~\ref{lem:measure_shrinking_lemma} to satisfy
\begin{align*}
     C \frac{\sigma}{\theta \lambda} < \nu, \quad \text{i.e.} \quad \sigma \le \frac{\nu \theta \lambda}{C}.
\end{align*}
This choice is admissible, since $C$ from Lemma~\ref{lem:measure_shrinking_lemma} is independent of $\sigma$. Furthermore, letting $\widetilde{\boldsymbol{\gamma}}$ be chosen in  Lemma~\ref{lem:a_de_giorgi_type_lemma},
we enforce that
\begin{align*}
   \max \{ \widetilde{\boldsymbol\gamma}, 1/\theta   \} \int_{t_1}^{t_2} \Vert f(\cdot,t) \Vert_{L^{\infty}(B_{8 \rho })} \,\dt \le  \sigma \epsilon \xi \boldsymbol{\om}, 
\end{align*}
and apply  Lemma~\ref{lem:measure_shrinking_lemma} with the above chosen $\sigma$, followed by Lemma~\ref{lem:a_de_giorgi_type_lemma} with $\xi$ replaced by $\sigma \epsilon \xi $, in the cylinder \begin{align*}
    (0, \overline{t}) + Q_{4 \rho}  ( \tfrac{1}{2}\theta   ) = B_{4 \rho} \times  ( \overline{t} - \tfrac{1}{2} \theta (4 \rho)^{2s}, \overline{t}  ]
\end{align*}
for an arbitrary $\overline{t} \in ( \frac{1}{2} \theta ( 4 \rho)^{2s}, \theta (4 \rho)^{2s}]$. Due to the arbitrariness of $\overline{t}$ we conclude that
\begin{align*}
    v \ge  \tfrac{1}{4} \sigma \epsilon \xi \boldsymbol{\om} \quad \text{a.e. in}\>\> B_{2 \rho} \times  ( \tfrac{1}{2} \theta (4 \rho )^{2s}, \theta ( 4 \rho)^{2s}  ].
\end{align*}
The proof is completed by setting $\eta = \frac{1}{4} \sigma \epsilon$, properly adjusting the relevant constants in dependence of the $\mathfrak{data}$, as well as $\lambda$ and relabeling $ \max\{ \widetilde{\boldsymbol\gamma}, 1/\theta \}$.

\section{Continuity of the gradient}
Recall that $v\equiv v_h$ is defined  in~\eqref{def:abbrevaition_for_diff_quotient}.
Before proving the continuity of $Du$, we exploit
Proposition~\ref{prop:expansion_of_positivity}  to prove a modulus of continuity of $v_h$ in the next result.

\begin{proposition}\label{prop:mod_of_cont_finite_difference_quotient}
    Let $u$ be a local weak solution to~\eqref{eq:PDE} 
    in the sense of Definition~\ref{def:weak_solution}.  
    Then, there exist $C>1$ and $\beta \in (0,1)$ depending only on the $\mathfrak{data}$, such that for $0<r<R$,  
    any $h \in \mathbb{R}^{N} \setminus \{0\}$ with $0 < |h| \le \frac{1}{32}R$, we have
    \begin{align*}
        \essosc_{Q_r(z_o) } v_h \le  C \bigg[  \boldsymbol{\om} \Big( \frac{r}{R} \Big)^{\beta} +  
        \int_{t_o - (rR)^{s}}^{t_o} \Vert f(\cdot,t) \Vert_{L^{\infty}(B_{\frac89 R}(x_o))}
        \, \dt \bigg].
    \end{align*}
    for some constant $C=C(\mathfrak{data})$ independent of $h$, where $f$ is defined  in~\eqref{def:widetilde_f_as_nonlocal_term}, i.e.
    \[
    f(x,t):=  \int_{\mathbb{R}^N \setminus B_{R}(x_o)} \frac{  v_h(y,t)}{\vert x-y \vert^{N+2s}}  
         \,\dy
    \]
    and
    \begin{align}\label{def:intermediate_definition_of_omega}
    \boldsymbol{\om} := \frac{1}{R} \int_{t_o -  (64R)^{2s}}^{t_o} \int_{\mathbb{R}^{N}} \frac{| u(x,t) |}{( R + | x - x_o |)^{N+2s}} \,\dx\dt.
\end{align}
\end{proposition}

\subsection{Proof of Proposition~\ref{prop:mod_of_cont_finite_difference_quotient}}

\subsubsection{The first step}\label{sec:the_first_step}
For $R>0$ let us introduce the cylinders at $z_o=(x_o, t_o)$:
\begin{align*}
    Q_o := 
    Q_{R}(z_o) \subset Q_{2R}(z_o) \subset
    \Omega_{T}.
\end{align*} 
For simplicity assume  $(x_o, t_o) =(0,0)$. Consider a quantity $\boldsymbol{\om}$ satisfying
\begin{align}\label{def:omega}
    \boldsymbol{\om} \ge 2 \esssup_{Q_{2 R} } | v |.
\end{align}
In the end, we will find that this $\boldsymbol{\om}$ can be selected as in \eqref{def:intermediate_definition_of_omega} apart from a multiplicative constant.
To continue, we define
\begin{align*}
    \boldsymbol{\mu}^{+} := \esssup_{Q_o} v, \quad \boldsymbol{\mu}^{-} := \essinf_{Q_o} v.
\end{align*}
We have from~\eqref{def:omega} the initial ocsillation control
\begin{align}\label{ineq:deriv_loc_bound_est_biggest_cyl}
    \essosc_{Q_o} v \le  \boldsymbol{\om}.
\end{align}
We will reduce the oscillation in a smaller cylinder contained in $Q_o$.
To this end, fix $\theta=\theta(\mathfrak{data}) \in (0,1)$ according to Proposition~\ref{prop:expansion_of_positivity} with $\lambda = \frac{1}{2}$; 
define 
\begin{align*}
    \tau := \theta  ( \tfrac19 R  )^{2s};
\end{align*}
consider the two alternatives:
\begin{align}\label{two_alternatives_for_mod_of_cont_proof}
    \left\{
	\begin{array}{c}
       \big| \big\{ v(\cdot, - \tau) - \boldsymbol{\mu}^{-} > \frac{1}{4}  \boldsymbol{\om} \big\} \cap B_{\frac19 R } \big| \ge \frac{1}{2} | B_{\frac19 R } |,
        \\ [5pt]
        \big| \big\{ \boldsymbol{\mu}^{+} - v(\cdot, - \tau)  > \frac{1}{4}  \boldsymbol{\om} \big\} \cap B_{\frac19 R  } \big| \ge \frac{1}{2} | B_{\frac19 R  } | .
        \end{array}\right.
\end{align}
We may assume $\boldsymbol{\mu}^{+} - \boldsymbol{\mu}^{-} \ge \frac{1}{2}  \boldsymbol{\om}$, then  one of the two in~\eqref{two_alternatives_for_mod_of_cont_proof} must hold, whereas the reverse case $\boldsymbol{\mu}^{+} - \boldsymbol{\mu}^{-} < \frac{1}{2}  \boldsymbol{\om}$ is favorable and automatically yields the aimed oscillation estimate \eqref{ineq:long_max_deriv_loc_bound_est_smaller_cyl}. Now, suppose that the first alternative holds; apply Proposition~\ref{prop:expansion_of_positivity} with $\lambda = \frac{1}{2}$, $\xi = \frac{1}{4}$ and $\rho = \frac19 R  $ to determine $\eta =\eta(\mathfrak{data})\in (0,1)$, as well as $t_1 = -(\frac19 R  )^{2s}$ and $t_2 = 0$, such that either  
\begin{align*}
    \widetilde{\boldsymbol\gamma} \int_{- (\frac19 R  )^{2s}}^{0} \Vert f(\cdot,t) \Vert_{L^{\infty}(B_{ \frac89 R  })} \,\dt >  \tfrac{1}{4}\eta \boldsymbol{\om}
\end{align*}
or
\begin{align*}
    v - \boldsymbol{\mu}^{-} \ge \tfrac14 \eta \boldsymbol{\om} \quad \text{a.e. in } Q_{\frac19 R}  ( \tfrac12 \theta ).
\end{align*}
Denote $\widehat{\boldsymbol\gamma} := 4 \widetilde{\boldsymbol\gamma}/\eta $ to simplify notation. Then, the above conclusion plus \eqref{ineq:deriv_loc_bound_est_biggest_cyl} yields that
\begin{align}\label{ineq:long_max_deriv_loc_bound_est_smaller_cyl}
    \essosc_{Q_{\frac19 R}(\frac12 \theta)} v &\le \max \bigg\{   ( 1-\tfrac{1}{4 }\eta  ) \boldsymbol{\om}, \widehat{\boldsymbol\gamma} \int_{ - (\frac19 R  )^{2s}}^{0} \Vert f(\cdot,t) \Vert_{L^{\infty}(B_{\frac89 R })} \,\dt 
    \bigg\} 
    =: \boldsymbol{\om}_{1}. 
\end{align}
Next, we introduce $R_{1} = \kappa R$ for some $\kappa \le \frac19$
in order to obtain the set inclusion  
\begin{align*}
    Q_{R_1} \subset Q_{ \frac19 R } ( \tfrac{1}{2}\theta  ), \quad \text{i.e.} \quad \kappa := \tfrac19 ( \tfrac{1}{2} \theta)^{\frac{1}{2}} .
\end{align*}
As a result of this set inclusion and~\eqref{ineq:long_max_deriv_loc_bound_est_smaller_cyl} we deduce
\begin{align*}
    Q_{R_1} \subset Q_o \quad \text{ and } \quad  \essosc_{Q_{R_1}} v \le  \boldsymbol{\om}_1,
\end{align*}
where the latter will play the role of~\eqref{ineq:deriv_loc_bound_est_biggest_cyl} in the next step. 
 Moreover, notice that the second alternative in~\eqref{two_alternatives_for_mod_of_cont_proof} can be dealt with analogously.

\subsubsection{The induction step}

Assume for $i =1,\dots,j$ that we have built
\begin{align*}
	\left\{
	\begin{array}{c}
R_o = R,\quad  R_i = \kappa R_{i-1},\quad \boldsymbol{\om}_o = \boldsymbol{\om}, 
\\[5pt]
B_i = B_{R_i}, \quad Q_i =Q_{R_i}, \quad \widetilde{Q}_{i} = B_{i} \times ( - (\frac19 R_i )^{2s}, 0], \quad Q_i \subset Q_{i-1} ,\\[5pt]
\dsty\boldsymbol{\mu}_{i}^{+} := \esssup_{Q_i} v, \quad \boldsymbol{\mu}_{i}^{-} := \essinf_{Q_i} v, \quad \essosc_{Q_{i}} v \le  \boldsymbol{\om}_i \\[5pt] \dsty \boldsymbol{\om}_{i} =  \max \bigg\{  (1-\tfrac{1}{4} \eta)  \boldsymbol{\om}_{i-1}, \widehat{\boldsymbol\gamma} \int_{ - ( \frac19 R_{i-1} )^{2s}}^{0} \Vert f(\cdot,t) \Vert_{L^{\infty}(B_{\frac89 R_{i-1} })} \,\dt 
\bigg\}.
\end{array}
	\right.
\end{align*}
Now we will show that the oscillation estimate continues to hold for the $(j+1)$-th step. 

Consider $\theta$ fixed as above. 
We repeat the argument from Section~\ref{sec:the_first_step} with $\boldsymbol{\mu}_{j}^{\pm}, \boldsymbol{\om}_{j}, Q_{j}$ and so on. Indeed, we define
\begin{align*}
    \tau_{j} := \theta ( \tfrac19 R_{j} )^{2s}
\end{align*}
and consider the two alternatives:
\begin{align}\label{iterated_two_alternatives_for_mod_of_cont_proof}
    \left\{
	\begin{array}{c}
       \big| \big\{ v(\cdot, - \tau_j)  - \boldsymbol{\mu}_{j}^{-} > \frac{1}{4}  \boldsymbol{\om}_{j} \big\} \cap B_{\frac19 R_j } \big| \ge \frac{1}{2} | B_{\frac19 R_j} |,
        \\ [5pt]
        \big| \big\{ \boldsymbol{\mu}_{j}^{+} - v(\cdot, - \tau_j)  > \frac{1}{4}  \boldsymbol{\om}_{j} \big\} \cap B_{\frac19 R_j} \big| \ge \frac{1}{2} | B_{\frac19 R_j} |. 
        \end{array}\right.
\end{align}
Like in the first step suppose that $\boldsymbol{\mu}_{j}^{+} - \boldsymbol{\mu}_{j}^{-} \ge \frac12 \boldsymbol{\om}_{j}$, such that one of the two alternatives is fulfilled. The case $\boldsymbol{\mu}_{j}^{+} - \boldsymbol{\mu}_{j}^{-} < \frac12  \boldsymbol{\om}_{j}$ can be incorporated into the forthcoming oscillation estimate~\eqref{ineq:inductive_long_max_deriv_loc_bound_est_smaller_cyl}. Without loss of generality assume that the first alternative in~\eqref{iterated_two_alternatives_for_mod_of_cont_proof} holds. We apply Proposition~\ref{prop:expansion_of_positivity} in $Q_j$, with $\lambda = \frac{1}{2}, \xi = \frac{1}{4}$ and $\rho = \frac19 R_j $ to get the same $\eta$ as in the first step, as well as $t_1 = - ( \frac19 R_j  )^{2s}$ and $t_2 = 0$, such that either 
\begin{align*}
    \widetilde{\boldsymbol\gamma} \int_{ - ( \frac19 R_j )^{2s}}^{0} \Vert f(\cdot,t) \Vert_{L^{\infty}(B_{\frac89 R_j })} \,\dt > \tfrac{1}{4}\eta \boldsymbol{\om}_{j}
\end{align*}
or 
\begin{align*}
    v - \boldsymbol{\mu}_{j}^{-} \ge \tfrac14 \eta \boldsymbol{\om}_{j} \quad \text{a.e. in}\>\> Q_{\frac19 R_j }  ( \tfrac12 \theta  ).
\end{align*}
Then, due to the $j$-th induction assumption we obtain that
\begin{align}\label{ineq:inductive_long_max_deriv_loc_bound_est_smaller_cyl}
    \essosc_{Q_{\frac19 R_j }(\frac12 \theta)} v &\le \max \bigg\{  (1-\tfrac{1}{4}\eta )  \boldsymbol{\om}_j
    , \widehat{\boldsymbol\gamma} \int_{- ( \frac19 R_j )^{2s}}^{0} \Vert f(\cdot,t) \Vert_{L^{\infty}(B_{\frac89 R_j} )} \,\dt 
    \bigg\} 
    =: \boldsymbol{\om}_{j+1}, 
\end{align}
where $\widehat{\boldsymbol\gamma}=\widehat{\boldsymbol\gamma}(\mathfrak{data})$ is as in~\eqref{ineq:long_max_deriv_loc_bound_est_smaller_cyl}. 
Now, let $R_{j+1} = \kappa R_{j}$ for the same $\kappa= \frac19  (\frac12 \theta)^{\frac12}$ chosen in the first step, in order
to get that
\begin{align}\label{condition_on_kappa_1}
     Q_{R_{j+1}} \subset Q_{\frac19 R_j }  ( \tfrac{1}{2}\theta  ).
\end{align}
As a consequence,~\eqref{condition_on_kappa_1} and~\eqref{ineq:inductive_long_max_deriv_loc_bound_est_smaller_cyl} yield
\begin{align*}
    Q_{R_{j+1}} \subset Q_{R_{j}} \quad \text{ and } \quad \essosc_{Q_{R_{j+1}}} v \le  \boldsymbol{\om}_{j+1},
\end{align*}
which finishes the induction step for the first alternative in~\eqref{iterated_two_alternatives_for_mod_of_cont_proof}. The same conclusion can be reached in case the second alternative of \eqref{iterated_two_alternatives_for_mod_of_cont_proof} holds true.

\subsubsection{Modulus of continuity}
From the previous section we have built  $R_n = \kappa^{n} R$, $Q_n = Q_{R_n}$ for $n \in \mathbb{N}$ such that $Q_n \subset Q_{n-1}$ and acquired the oscillation estimate
\begin{align}\label{ineq:n_th_step__long_max_deriv_loc_bound_est_smaller_cyl}
    \essosc_{Q_n } v &\le \max \bigg\{ (1-\tfrac14 \eta ) \boldsymbol{\om}_{n-1} 
    , \widehat{\boldsymbol\gamma}  \int_{ - ( \frac19 R_{n-1} )^{2s}}^{0} \Vert f(\cdot,t) \Vert_{L^{\infty}(B_{\frac89 R_{n-1} })} \,\dt 
    \bigg\} \\ \nonumber
    &=: \boldsymbol{\om}_{n}. 
\end{align}
Thus, 
by the definition of $\boldsymbol{\om}_n$ we obtain the recursive estimate
\begin{align*}
    \boldsymbol{\om}_n &=
    \max\bigg\{ 
    (1-\tfrac{1}{4}\eta )  \boldsymbol{\om}_{n-1}, C \int_{ - ( \frac19 R_{n-1} )^{2s}}^{0} \Vert f(\cdot,t) \Vert_{L^{\infty}(B_{\frac89 R_{n-1} })} \,\dt \bigg\} \\
    &\le 
    \max\bigg\{ 
    (1-\tfrac{1}{4}\eta )  \boldsymbol{\om}_{n-1}, C \int_{ - (\frac19 R )^{2s}}^{0} \Vert f(\cdot,t) \Vert_{L^{\infty}(B_{\frac89 R })} \,\dt \bigg\}.
\end{align*}
Iterating the above inequality, while considering the estimate in~\eqref{ineq:n_th_step__long_max_deriv_loc_bound_est_smaller_cyl} yields
\begin{align}\label{ineq:completed_iteration_mod_od_cont_proof_no_radii}
    \essosc_{Q_n } v \le \boldsymbol{\om}_n \le (1-\tfrac{1}{4}\eta )^{n}  \boldsymbol{\om} + C
    \int_{ - R^{2s}}^{0} \Vert f(\cdot,t) \Vert_{L^{\infty}(B_{\frac89 R })}
   \, \dt.
\end{align}
Since  $\{ R_{n} \}_{n \in \mathbb{N}_0 }$ partitions $(0,R)$, for $r \in (0,R)$, there exists $n \in \mathbb{N}_0$ such that
\begin{align}\label{introduction_of_r_in_mod_of_continuity}
    R_{n+1} \le r < R_n.
\end{align}
The second inequality in \eqref{introduction_of_r_in_mod_of_continuity} implies that
\begin{align*}
    Q_r \subset Q_n .
\end{align*}
Using the definition of $R_n$ and the first inequality in~\eqref{introduction_of_r_in_mod_of_continuity} we obtain that
\begin{align*}
    r \ge R_{n+1} = \kappa^{n+1} R,
\end{align*}
which gives us that
\begin{align}\label{radii_iteration_bound_for_first_eta_term}
    \frac{r}{R} \ge  \kappa^{n+1} \quad \implies \quad (1-\tfrac14 \eta )^{n+1} \le \Big( \frac{r}{R} \Big)^{\beta},
\end{align}
where
\begin{align*}
    \beta = \frac{\ln(1- \frac14 \eta)}{\ln \kappa }.
\end{align*}
By inserting~\eqref{radii_iteration_bound_for_first_eta_term} into~\eqref{ineq:completed_iteration_mod_od_cont_proof_no_radii} we obtain for any $r \in (0,R)$ that
\begin{align}\label{final_in_proof_oscillation_estimate}
    \essosc_{Q_r } v \le  2  \boldsymbol{\om} \Big( \frac{r}{R} \Big)^{\beta} + C
    \int_{ - R^{2s}}^{0}  \Vert f(\cdot,t) \Vert_{L^{\infty}(B_{\frac89 R  })}
    \,\dt.
\end{align}
Now, we perform an interpolation argument. Indeed, we may assume that~\eqref{final_in_proof_oscillation_estimate} holds true, instead of $R$, for some $\overline{R} \in (r, R)$, and then setting $\overline{R} = ( r R) ^{\frac{1}{2}}$ to get that 
\begin{align}\label{ineq:penultimate_ess_osc_bound_for_finite_diff_quotient_in_mod_cont_proof}
    \essosc_{Q_r } v \le  2  \boldsymbol{\om} \Big( \frac{r}{R} \Big)^{\frac{\beta}{2}} + C
    \int_{ - (rR)^{s}}^{0} \Vert f(\cdot,t) \Vert_{L^{\infty}(B_{\frac89 R})}
   \, \dt.
\end{align}
This conclusion is reached by assuming $\boldsymbol{\om}$ satisfies \eqref{def:omega}.
Now, we choose $\boldsymbol{\om}$ to satisfy \eqref{def:omega}. In fact, using \eqref{ineq:abs_value_of_finite_difference_quotient_bound} with $y = x+h$, while dividing by $|h|$, as well as taking the supremum over $Q_{2R}$,
and then \cite[Theorem 1.6]{Liao_24_1} with $b=0$, we observe that
\begin{align*}
    \esssup_{Q_{2 R}} | v | &\le \frac{C}{R} \Big[ \Vert u \Vert_{L^{\infty}(Q_{16 R})} + \mathrm{Tail}(u; Q_{16 R} ) \Big] \\
    &\le \frac{C}{R} \int_{ - (64 R)^{2s}}^{0} \int_{\mathbb{R}^{N}} \frac{| u(x,t) |}{( R + | x |)^{N+2s}} \,\dx\dt.
\end{align*}  
Hence, to fulfill \eqref{def:omega} we only need to pick 
\begin{align*}
    \boldsymbol{\om} = \frac{2C}{R} \int_{ -  (64R)^{2s}}^{0} \int_{\mathbb{R}^{N}} \frac{| u(x,t) |}{( R + | x |)^{N+2s}} \,\dx\dt,
\end{align*}
and then \eqref{ineq:penultimate_ess_osc_bound_for_finite_diff_quotient_in_mod_cont_proof} holds true.
Finally, we substitute the above $\boldsymbol{\om}$ in \eqref{ineq:penultimate_ess_osc_bound_for_finite_diff_quotient_in_mod_cont_proof}, redefine $\frac{\beta}{2}$ to be $\beta$, properly enlarge $C$,
and conclude the proof of Proposition~\ref{prop:mod_of_cont_finite_difference_quotient} with the $\boldsymbol{\om}$ in \eqref{def:intermediate_definition_of_omega}.

\subsection{Proof of Theorem~\ref{Thm:main}}

We conclude similar as in~\cite[Proof of Theorem 1.6]{DKLN_2025}.
First, choose $h = \overline{h} e_{i}$, where $\overline{h} \in (0, \frac{1}{32} R)$ and $\{e_i:\, i=1,2,\dots,N\}$ are canonical unit vectors in $\mathbb{R}^{N}$.
    By Proposition~\ref{prop:mod_of_cont_finite_difference_quotient} 
    we obtain for $r \in (0,R)$ small enough that
    \begin{align*}
        \essosc_{Q_r (x_o, t_o)} v_{h} \le  C    
        \bigg[  \boldsymbol{\om} \Big( \frac{r}{R} \Big)^{\beta} +  \int_{t_o - (rR)^{s}}^{t_o} \Vert f(\cdot,t) \Vert_{L^{\infty}(B_{\frac89  R}(x_o))}        
        \,\dt \bigg],
    \end{align*}
    for some constant $C=C(\mathfrak{data})$. 
    Next, by applying Lemma~\ref{lem:first_tail_estimate} for the integrand of the last integral, followed by \cite[Theorem 1.6]{Liao_24_1} with $b=0$, we deduce 
    \begin{align}\label{essosc_of_sequence_v_h}
        \essosc_{Q_r (x_o, t_o) } v_{h} &\le  C    
        \bigg[  \boldsymbol{\om} \Big( \frac{r}{R} \Big)^{\beta} + \frac{1}{R^{2s+1}} \int_{t_o - (rR)^{s}}^{t_o} \mathcal{T}(t; B_{2R}(x_o) )      
        \,\dt \bigg] \nonumber \\
        &\le C \bigg[  \boldsymbol{\om} \Big( \frac{r}{R} \Big)^{\beta} + \frac{1}{R} \int_{t_o - (rR)^{s}}^{t_o} \int_{\mathbb{R}^{N}} \frac{| u(x,t) |}{( R + | x - x_o |)^{N+2s}} \,\dx  \,\dt \bigg].
    \end{align}
    Note that the right-hand side does not depend on $h$ anymore. 
    Now, consider the sequence $\{ v_{h_n} \}_{n \in \mathbb{N}}$, where $h_{n} \to 0$ as $n \to \infty$, satisfying~\eqref{essosc_of_sequence_v_h}. First of all, notice that 
    this sequence is uniformly bounded, 
    since for any $z \in \overline{Q}_{\frac12 r} (x_o, t_o) := \overline{B}_{\frac12 r} (x_o) \times [t_o - ( \frac{1}{2} r)^{2s}, t_o]$ we have
    \begin{align}\label{ineq:unif_bd}
        | v_{h_n} (z) | \le  \Vert Du \Vert_{L^{\infty}(Q_{2 R})} \le \frac{C}{R} \Big[\|u\|_{L^\infty(Q_{16 R} )} +  \mathrm{Tail}(u;Q_{16 R} ) \Big],
    \end{align}
    where for the last inequality we used Proposition~\ref{prop:Lipschitz:1}. 
    Observe that the right-hand side is independent of $z$ and $h_n$.
    Moreover,  
    the sequence $\{ v_{h_n} \}_{n \in \mathbb{N}}$ is equicontinuous, 
    since for any $z_1, z_2 \in \overline{Q}_{\frac12 r} (x_o, t_o)$ we have
    \begin{align}\label{ineq:equicont}
        | v_{h_n} (z_1) - v_{h_n} (z_2) | &\le  \essosc_{Q_r (x_o, t_o) } v_{h_n} \nonumber \\
        &\le  C \bigg[  \boldsymbol{\om} \Big( \frac{r}{R} \Big)^{\beta} + \frac{1}{R} \int_{t_o - (rR)^{s}}^{t_o} \int_{\mathbb{R}^{N}} \frac{| u(x,t) |}{( R + | x - x_o |)^{N+2s}} \,\dx  \,\dt \bigg],
    \end{align}
    due to~\eqref{essosc_of_sequence_v_h}. 
    Combining~\eqref{ineq:unif_bd} and~\eqref{ineq:equicont} enables us to apply Arzel\`a-Ascoli's theorem in order to obtain a subsequence $\{ v_{h_{n_k}} \}_{k \in \mathbb{N}}$, which converges uniformly to $D_i u$. 
    Then, taking the limit $h_{n_k} \to 0$ in~\eqref{essosc_of_sequence_v_h} on the smaller domain $Q_{\frac{1}{2} r}$ gives us the estimate
    \begin{align*}
        \essosc_{Q_{\frac{1}{2} r} (x_o, t_o)} D_{i} u \le  C \bigg[  \boldsymbol{\om} \Big( \frac{r}{R} \Big)^{\beta} + \frac{1}{R} \int_{t_o - (rR)^{s}}^{t_o} \int_{\mathbb{R}^{N}} \frac{| u(x,t) |}{( R + | x - x_o  |)^{N+2s}} \,\dx\dt \bigg],
    \end{align*}
for some constant $C=C(\mathfrak{data})$. 
A scaling and covering argument completes the proof with $\boldsymbol\omega$ as in~\eqref{def:final_omega}.

\

\noindent{\bf Data Availability.} Data sharing not applicable to this article as no datasets were generated or analyzed during the current study.

\

\noindent{\bf Conflict of Interest.} The author declares no conflict of interest.



\begin{thebibliography}{99}


\bibitem{BDLGMS24}
\newblock V.~B\"ogelein, F.~Duzaar, N.~Liao, G.~Molica Bisci and R.~Servadei.
\newblock Regularity for the fractional $ p $-Laplace equation.
\newblock \emph{J. Funct. Anal.}, {\bf289}(9), (2025),  Paper No. 111078, 69 pp.  


\bibitem{BDLGM24-2}
\newblock V.~B\"ogelein, F.~Duzaar, N.~Liao, G.~Molica Bisci and R.~Servadei.
\newblock Gradient regularity for $(s,p)$-harmonic functions. 
\newblock \emph{Calc. Var. Partial Differential Equations}, {\bf 64}(8), (2025), Paper No. 253, 63 pp.


\bibitem{BL-17}
\newblock L.~Brasco and E.~Lindgren.
\newblock Higher Sobolev regularity for the fractional $p$-Laplace equation in the superquadratic case.
\newblock \emph{Adv. Math.}, {\bf 304}, (2017), 300--354.


\bibitem{BLS-18}
\newblock L.~Brasco, E.~Lindgren and A.~Schikorra.
\newblock Higher H\"older regularity for the fractional $p$-Laplacian in the superquadratic case.
\newblock \emph{Adv. Math.}, {\bf 338}, (2018), 782--846.


\bibitem{CCV_2011}
\newblock L.~A.~Caffarelli, C.~H.~Chan and A.~F.~Vasseur.
\newblock Regularity theory for parabolic nonlinear integral operators.
\newblock \emph{J. Amer. Math. Soc.}, {\bf 24}(3), (2011) 849--869.


\bibitem{CLD_2016}
\newblock H.~A.~Chang-Lara and G.~D\'avila.
\newblock Hölder estimates for non-local parabolic equations with critical drift. 
 \newblock \emph{J. Differential Equations}, {\bf 260}(5), (2016), 4237--4284.


\bibitem{DeFilippis}
\newblock F.~De Filippis.
\newblock Local boundedness for weak solutions to fractional porous medium equation.
\newblock arXiv:2602.23001.  


\bibitem{DiBenedetto}
\newblock E.~DiBenedetto.
\newblock Degenerate parabolic equations.
\newblock \emph{Universitext}, Springer, New York, 1993.


\bibitem{DB-RA}
\newblock E.~DiBenedetto. 
\newblock Real analysis. 2nd Ed., \emph{Birkhäuser Advanced Texts}, Birkhäuser/Springer, New York, 2016. 


\bibitem{DKLN_2025}
\newblock L.~Diening, K.~Kim, H.S.~Lee, and S.~Nowak.
\newblock Gradient estimates for parabolic nonlinear nonlocal equations.
\newblock \emph{Calc. Var. Partial Differential Equations}, {\bf 64}(3),  (2025), Paper No. 98, 86 pp.


\bibitem{Diening_Kim_Lee_Nowak_24}
\newblock L.~Diening, K.~Kim, H.S.~Lee, and S.~Nowak.
\newblock Nonlinear nonlocal potential theory at the gradient level.
\newblock arXiv preprint arXiv:2402.04809.


\bibitem{DKLN_ann}
\newblock L.~Diening, K.~Kim, H.S.~Lee, and S.~Nowak. 
\newblock Higher differentiability for the fractional $p$-Laplacian. 
\newblock \emph{Math. Ann.}, {\bf 391}(4), (2025), 5631--5693.


\bibitem{Domokos_2004_2}
\newblock A.~Domokos.
\newblock Differentiability of solutions for the non-degenerate $p$-Laplacian in the Heisenberg group.
\newblock \emph{J. Differential Equations}, {\bf 204}(2), (2024), 439--470.


\bibitem{KW_24}
\newblock M.~Kassmann and M.~Weidner.
\newblock The parabolic Harnack inequality for nonlocal equations.
\newblock \emph{Duke Math. J.}, {\bf 173}(17), (2024), 3413--3451.


\bibitem{LSU}
\newblock O. A. Lady\v{z}enskaja, V. A. Solonnikov, and N. N. Ural'ceva.
\newblock \emph{Linear and quasilinear equations of parabolic type.} (Russian) Translated from the Russian by S. Smith. Translations of Mathematical Monographs, Vol. 23. American Mathematical Society, Providence, RI, 1968. xi+648 pp.


\bibitem{Liao_24_1}
\newblock N.~Liao.
\newblock Harnack estimates for nonlocal drift-diffusion equations.
\newblock arXiv:2402.11986v2.


\bibitem{Liao_24_2}
\newblock N.~Liao.
\newblock H\"older regularity for the parabolic fractional $p$-Laplacian.
\newblock \emph{ Calc. Var. Partial Differential Equations} {\bf 63}(1), (2024), Paper No. 22, 34 pp.


\bibitem{Liao_24_3}
\newblock N.~Liao.
\newblock On the modulus of continuity of solutions to nonlocal parabolic equations.
\newblock \emph{J. Lond. Math. Soc. (2)} {\bf 110}(3), (2024), Paper No. e12985, 30 pp.


\bibitem{Stein}
\newblock E.~M. Stein.
\newblock Singular integrals and differentiability properties of functions.
\newblock \emph{Princeton Math. Ser.} {\bf 30}, Princeton Univ. Press, Princeton, NJ, (1970).


\end{thebibliography}
\end{document}